\documentclass[tbtags,reqno]{amsart}

\edef\savecatcodeat{\the\catcode`@} \catcode`\@=11

\def\tb@ifSpecChars#1#2{#1}
\def\tb@ifNoSpecChars#1#2{#2}

\def\tableau{%
  \bgroup% matched in \tb@tableauD
  \@ifstar{\let\Tif\tb@ifNoSpecChars\tb@tableauB}% *, don't use special chars
          {\let\Tif\tb@ifSpecChars\tb@tableauB}}% no *, use special chars

\def\tb@tableauB{% add [] if no [options]
  \@ifnextchar[{\tb@tableauC}{\tb@tableauC[]}}

\def\tb@tableauC[#1]{\hbox\bgroup%
    \let\\=\cr% end line
    \def\bl{\global\let\tbcellF\tb@cellNF}%
    \def\tf{\global\let\tbcellF\tb@cellH}% highlighted cell
    \dimen2=\ht\strutbox \advance\dimen2 by\dp\strutbox%
    \ifx\baselinestretch\undefined\relax%
    \else%
       \dimen0=100sp \dimen0=\baselinestretch\dimen0%
       \dimen2=100\dimen2 \divide\dimen2 by\dimen0%
    \fi%
    \let\tpos\tb@vcenter% default position
    \tb@initYoung% default tableau type
    \tb@options#1\eoo% parse options
    \let\arrow\tb@arrow%
    \dimen0=\Tscale\dimen2%
    \dimen1=\dimen0 \advance\dimen1 by \tb@fframe%
    \lineskip=0pt\baselineskip=0pt% line spacing will be from \vbox to \dimen0
    \def\tb@nothing{}%
    \def\endcellno{$\rss\egroup\bss\egroup}% end cell w/o overlap
    \def\endcell{\endcellno\kern-\dimen0}% end cell & prepare to overlap it
    \def\begincell{\vbox to\dimen0\bgroup\vss\hbox to\dimen0\bgroup\hss$}%
    \let\overlay\tb@overlay%
    \let\fl\tb@fl%
    \let\lss\hss\let\rss\hss\let\tss\vss\let\bss\vss% cell alignment
    \def\mkcell##1{% format individual cell
        \let\tbcellF\tb@cellD% default cell frame
        \def\tb@cellarg{##1}% store cell contents
        \ifx\tb@cellarg\tb@nothing\let\tb@cellarg\tb@cellE\fi%
        \begincell\tb@cellarg\endcellno% the actual cell content
        \tbcellF}% draw cell frame
    \let\savecellF\tbcellF% save global value of cellF in case of nested tableau
     \Tif{\catcode`,=4\catcode`|=\active}{}\tb@tableauD}%

\let\tb@savetableauD\tableauD% save any current definition
{% set up characters which will be interpreted as command characters
    \catcode`|=\active \catcode`*=\active \catcode`~=\active%
    \catcode`@=\active% command characters
\gdef\tableauD#1{%
  \Tif{% make all the command characters active in math mode when #1 parsed
    \mathcode`|="8000 \mathcode`*="8000%
    \mathcode`~="8000 \mathcode`@="8000%
    \def@{\bullet}%
    \let|\cr% end line
    \let*\tf% highlighted cell
    \let~\sk% skew cell
  }{}%
  \tpos{\tabskip=0pt\halign{&\mkcell{##}\cr#1\crcr}}%
  \global\let\tbcellF\savecellF% restore global value
  \egroup% match \hbox\bgroup at start of \tableauC
  \egroup}% match \bgroup at start of \tableau
}
\let\tb@tableauD\tableauD% rename the command
\let\tableauD\tb@savetableauD% restore old command with this name
\let\tb@savetableauD\undefined

\def\tb@options#1{\ifx#1\eoo\relax\else\tb@option#1\expandafter\tb@options\fi}

\def\tb@option#1{%
  \if#1t\let\tpos\tb@vtop\fi%        t = align at top
  \if#1c\let\tpos\tb@vcenter\fi%     c = align at center
  \if#1b\let\tpos\vbox\fi%           b = align at bottom
  \if#1F\tb@initFerrers\fi%          F = Ferrers diagram
  \if#1Y\tb@initYoung\fi%            Y = Young diagram
  \if#1s\tb@initSmall\fi%            s = small boxes
  \if#1m\tb@initMedium\fi%           m = medium boxes
  \if#1l\tb@initLarge\fi%            l = large boxes
  \if#1p\tb@initPartition\fi%            p = small partition sized boxes
  \if#1a\tb@initArrow\fi%            a = use arrow font as base dimension
}

\def\tb@vcenter#1{\ifmmode\vcenter{#1}\else$\vcenter{#1}$\fi}

\def\tb@vtop#1{\hbox{\raise\ht\strutbox\hbox{\lower\dimen0\vtop{#1}}}}

\def\tb@initPartition{\def\Tscale{.3}}
\def\tb@initSmall{\def\Tscale{1}}
\def\tb@initMedium{\def\Tscale{2}}
\def\tb@initLarge{\def\Tscale{3}}

\def\tb@initArrow{\dimen2=1.25em}

\def\tb@initYoung{%
  \def\tb@cellE{}% empty cell stays empty
  \let\tb@cellD\tb@cellN% default frame is normal frame
  \def\sk{\global\let\tbcellF\tb@cellNF}}% skew cells are empty
\def\tb@initFerrers{%
  \def\tb@cellE{\bullet}% empty cell gets bullet
  \let\tb@cellD\tb@cellNF% default frame is no frame
  \def\sk{\bullet}}% skew cell gets bullet

\tb@initMedium% default scale

\def\tb@sframe#1{%
  \vbox to0pt{%            Embed frame in a box of no vert or hor extent
    \vss%                            pull box above reference point
    \hbox to0pt{%
      \hss%                          pull box left of reference point
      \vbox to\dimen1{%              Actual width of frame
        \hrule depth #1 height0pt% draw top edge of frame
        \vss%                     begin vcenter sides
        \hbox to\dimen1{%           horiz box with side edges just inside
          \vrule width #1 height\dimen1% left edge
          \hss%                     stretch center
          \vrule width #1%         right edge
          }%
        \vss%                     end vcenter sides
        \hrule height #1 depth 0in% bottom edge
        }%
      \kern-\tb@hframe%           horiz alignment off by half line width
      }%
    \kern-\tb@hframe}}%           vert alignment off by half line width
\def\tb@hframe{.2pt}\def\tb@fframe{.4pt}\def\tb@bframe{2pt}
\def\tb@cellH{\tb@sframe{\tb@bframe}}       % bold frame
\def\tb@cellNF{}                            % no frame
\def\tb@cellN{\tb@sframe{\tb@fframe}}       % normal frame
\let\tbcellF\tb@cellN                       % default is normal

\def\tb@rpad{1pt}
\def\tb@lpad{1pt}
\def\tb@tpad{1.8pt}
\def\tb@bpad{1.8pt}

\def\tb@overlay{\endcell\@ifnextchar[{\tb@overlaya}{\begincell}}
\def\tb@overlaya[#1]{\vbox to\dimen0\bgroup%
  \tb@overlayoptions#1\eoo%
  \tss\hbox to\dimen0\bgroup\lss$}
\def\tb@overlayoptions#1{\ifx#1\eoo\relax\else\tb@overlayoption#1\expandafter\tb@overlayoptions\fi}

\def\tb@overlayoption#1{
  \if#1t\def\tss{\vskip\tb@tpad}\let\bss\vss\fi% t = align at top
  \if#1c\let\tss\vss\let\bss\vss\fi%             c = align at center
  \if#1b\def\bss{\vskip\tb@bpad}\let\tss\vss\fi% b = align at bottom
  \if#1l\def\lss{\hskip\tb@lpad}\let\rss\hss\fi% l = align at left
  \if#1m\let\lss\hss\let\rss\hss\fi%             m = align at middle
  \if#1r\def\rss{\hskip\tb@rpad}\let\lss\hss\fi% r = align at right
}

\def\tb@fl{\endcell\begincell\vrule depth 0pt width \dimen0 height \dimen0 \endcell\begincell}

\@ifundefined{diagram}{}{
\def\tb@arrowpad{.5}

\newoptcommand{\tb@arrow}{\@ne}[2]{%
  \endcell% end previous cell contents
   \begingroup%
   \let\dg@getnodesize\tb@getnodesize% substitute routine to get nodesize
   \dg@USERSIZE=#1\relax%
   \ifnum\dg@USERSIZE<\@ne \dg@USERSIZE=\@ne \fi%
   \dg@parse{#2}%
   \dg@label{\tb@draw{#1}{#2}}}% draw arrow

\def\tb@getnodesize#1#2#3#4#5{\dimen3=\tb@arrowpad\dimen2 #4=\dimen3 #5=\dimen3\relax}
\def\tb@getnodesize#1#2#3#4#5{\ifnum#2=0\ifnum#3=0\tb@getnodesizetail{#4}{#5}\else\tb@getnodesizehead{#4}{#5}\fi\else\tb@getnodesizehead{#4}{#5}\fi}
\def\tb@getnodesizetail#1#2{\dimen3=.5\dimen2 #1=\dimen3 #2=\dimen3}
\def\tb@getnodesizehead#1#2{\dimen3=.5\dimen2 #1=\dimen3 #2=\dimen3}

\def\tb@draw#1#2#3#4{%
        \dg@X=0\dg@Y=0\dg@XGRID=1\dg@YGRID=1\unitlength=.001\dimen0%
        \dg@LBLOFF=\dgLABELOFFSET \divide\dg@LBLOFF\unitlength%
        \dg@drawcalc% compute arrow geometry
        \begincell% start tableau cell
        \let\lams@arrow\tb@lams@arrow% substitute routine
        \begin{picture}(0,0)\begingroup\dg@draw{#1}{#2}{#3}{#4}\end{picture}%
        \endcell% end tableau cell
        \endgroup% match \begingroup in \tb@arrow
        \begincell}% start new entry in this cell
}

\def\tb@lams@arrow#1#2{%
 \lams@firstx\z@\lams@firsty\z@
 \lams@lastx#1\relax\lams@lasty#2\relax
 \lams@center\z@
 \N@false\E@false\H@false\V@false
 \ifdim\lams@lastx>\z@\E@true\fi
 \ifdim\lams@lastx=\z@\V@true\fi
 \ifdim\lams@lasty>\z@\N@true\fi
 \ifdim\lams@lasty=\z@\H@true\fi
 \NESW@false
 \ifN@\ifE@\NESW@true\fi\else\ifE@\else\NESW@true\fi\fi
 \ifH@\else\ifV@\else
  \lams@slope
  \ifnum\lams@tani>\lams@tanii
   \lams@ht\ten@\p@\lams@wd\ten@\p@
   \multiply\lams@wd\lams@tanii\divide\lams@wd\lams@tani
  \else
   \lams@wd\ten@\p@\lams@ht\ten@\p@
   \divide\lams@ht\lams@tanii\multiply\lams@ht\lams@tani
  \fi
 \fi\fi
 \ifH@  \lams@harrow
 \else\ifV@ \lams@varrow
 \else \lams@darrow
 \fi\fi
}

\catcode`\@=\savecatcodeat
\let\savecatcodeat\undefined

\usepackage{mathrsfs}
\usepackage{latexsym}
\usepackage{amsthm,amsopn,amsfonts,amssymb,epsfig}
\numberwithin{equation}{section}
\makeatletter
\renewcommand{\subsubsection}{\@startsection
{subsubsection} {3} {0mm} {\baselineskip} {-0.5\baselineskip}
{\normalfont\normalsize\bfseries}} \makeatother

\newtheorem{theorem}{Theorem}
\newtheorem{lemma}[theorem]{Lemma}
\newtheorem{proposition}[theorem]{Proposition}

\newtheorem{corollary}[theorem]{Corollary}

\newtheorem{definition}[theorem]{Definition}

\newtheorem{remark}[theorem]{Remark}
\newtheorem*{acknow}{Acknowledgments}

\def\la{{\lambda}}

\def\cal L{{\mathcal L}}

\newcommand{\cercle}[1]{\ensuremath{\setlength{\unitlength}{1ex}\begin{picture}(2.8,2.8)\put(1.4,0.7){\circle{2.8}\makebox(-5.6,0){#1}}\end{picture}}}
\newcommand{\tcercle}[1]{\ensuremath{\setlength{\unitlength}{1ex}\begin{picture}(2.8,2.8)\put(1.4,1.4){\circle{2.8}\makebox(-5.6,0){#1}}\end{picture}}}

\newcommand{\ds}{\displaystyle}
\newcommand{\mc}{\mathcal}

\newcommand{\spar}{\ensuremath{\mbox{\textnormal{SPar}}}}
\newcommand{\versg}[1]{\ensuremath{\overleftarrow{\,\phantom{|}#1\,\phantom{|}}}}
\newcommand{\versd}[1]{\ensuremath{\overrightarrow{\,\phantom{|}#1\,\phantom{|}}}}
\newcommand{\aand}{\ensuremath{ \quad\textrm{and}\quad}}
\newcommand{\wwhere}{\ensuremath{ \quad\textrm{where}\quad}}
\let\d\partial
\let\n\noindent

\def\M{{\mathcal M}}

\let\la\lambda
\let\La\Lambda
\let\Om\Omega
\let\om\omega
\let\ta\theta

\let\rw\rightarrow

\newcommand{\LL}{\ensuremath{\langle\!\langle}}
\newcommand{\RR}{\ensuremath{\rangle\!\rangle}}

\begin{document}

\title{Classical symmetric functions in superspace}

\author{Patrick Desrosiers} \thanks{P.Desrosiers@ms.unimelb.edu.au}
\address{Department of mathematics and statistics, The University of Melbourne, Parkville, Australia, 3010.}
\author{Luc Lapointe}
\thanks{lapointe@inst-mat.utalca.cl}
\address{Instituto de Matem\'atica y F\'{\i}sica, Universidad de
Talca, Casilla 747, Talca, Chile.}
\author{Pierre Mathieu} \thanks{pmathieu@phy.ulaval.ca}
\address{D\'epartement de physique, de g\'enie physique et
d'optique, Universit\'e Laval,  Qu\'ebec, Canada, G1K 7P4.}

 \maketitle

\begin{quote}
\small{\sc Abstract.} We present  the basic elements of a
generalization of symmetric function theory  involving functions  of
commuting and anticommuting (Grassmannian) variables. These new
functions, called symmetric functions in superspace, are invariant
under the diagonal action of the symmetric group on the sets of
commuting and anticommuting variables.  In this work, we present the
superspace extension of the classical bases, namely, the monomial
symmetric functions, the elementary symmetric functions, the
completely symmetric functions, and the power sums. Various basic
results, such as the generating functions for the multiplicative
bases, Cauchy formulas, involution operations as well as the
combinatorial scalar product are also generalized.
% Finally a one-parameter  extension of this scalar
% product as well as the completely symmetric functions are introduced.
\end{quote}

%\tableofcontents

\section{Introduction}

% \subsection{Symmetric superpolynomials}

Superspace is an extension of Euclidean space in $N$ variables
involving anticommuting variables. Its coordinates $(x_1, \cdots
x_N, \ta_1, \cdots \ta_N)$ obey the relations $x_ix_j=x_jx_i, \,
x_i\ta_j=\ta_jx_i$, and $\ta_i\ta_j=-\ta_j\ta_i$. Functions in
superspace, also called superfunctions,  are thus functions of two
types of variables. For instance, when $N=2$, all functions are
combinations of the following expressions
\begin{equation}
f_0(x_1,x_2)\, , \qquad \ta_1 f_1(x_1,x_2) + \ta_2 f_2(x_1,x_2)\, ,
\qquad \ta_1\ta_2 f_{3}(x_1,x_2)\, ,
\end{equation}
where the $f_i$'s stand for arbitrary functions of $x_1$ and $x_2$.
Functions of the second type are fermionic (alternatively said to be
odd) while those of the first and third types are bosonic (even).

The aim of this work is to lay down the foundation of a symmetric
function theory in superspace, where by a symmetric function in
superspace, we understand
%L is not a function invariant
%L under each type of symmetry transformation separately
%L (which would be too restrictive since  $\theta_1+\ldots+\theta_N$
%L is the only symmetric function of $\theta$), but rather
a function invariant  under the simultaneous interchange of $x_i
\leftrightarrow x_j$ and $\theta_i \leftrightarrow \theta_j$ for any
$i,j$. Examples when $N=2$ of symmetric polynomials in superspace
are
\begin{equation}
x_1^2x_2^2\, , \quad\ta_1 x_1^4+\ta_2x_2^4\, , \quad \ta_1x_2^2+
\ta_2x_1^2\,, \quad  \ta_1\ta_2(x_1^3x_2-x_1x_2^3)\,.
\end{equation}
The enforced interconnection between the transformation properties
of the bosonic and the fermionic variables is what makes the
resulting objects most interesting and novel.

The first step in the elaboration of a theory of symmetric
polynomials in superspace is the introduction of a proper labeling
for bases of the ring of symmetric superpolynomials. This
generalization of the concept of partitions, which was called
superpartition in \cite{DLM1}, turns out to be equivalent to what is
known as a MacMahon standard diagram in \cite{pak}. With this
concept in hand, the construction of the superextension of the
symmetric monomial basis (supermonomial basis for short) is rather
immediate \cite{DLM1}.
  From there on, the  natural route
for extending to superspace the multiplicative classical symmetric
functions, such as elementary, homogeneous and power sum symmetric
functions, is via the extension of their generating functions.
%L For the elementary $e_n$, homogeneous $h_n$ and power sum
%L $p_n$ symmetric functions, these generating functions are
%L respectively given by \cite{Mac}:
%L \begin{equation}
%L \sum_{n\geq 0} e_n\, t^n=\prod_{i\geq 1}(1+x_it)\; ,\quad
%L \sum_{n\geq 0} h_n\, t^n=\prod_{i\geq 1}\frac{1}{(1-x_it)}\;,\quad
%L \sum_{n\geq 1} p_n\, t^n=\prod_{i\geq 1}\frac{x_i t}{(1-x_it) }\, .
%L \end{equation}
%L The basis elements are generated from the product of these
%L functions:
%L \begin{equation}
%L f_\la = f_{\la_1}\cdots f_{\la_n}
%L \end{equation}
%L where $\la$ denotes a partition and $f$ is any of $e,\,h$ or $p$.
The central point of this extension lies in the observation that the
replacement $tx_i\rw tx_i+\tau \theta_i$,  where $t$ is the usual
counting variable and  $\tau$ is an anticommuting parameter, which
lifts  the generating functions directly  to superspace, yields the
``appropriate'' bases.  That is, the bases that are obtained have
properties that extend those satisfied by their classical
counterparts, such as for instance orthogonality relation and
determinantal formulas.  An even more convincing argument as to why
this is the right symmetric function theory in superspace comes from
its connection with an $N$-body problem in supersymmetric quantum
mechanics involving a parameter $\beta$ (see e.g., \cite{DLM1,DLM5}
and references therein).  The eigenfunctions of the Hamiltonian  of
this model are superspace generalizations of Jack polynomials that
specialize to various of the  bases presented in this article, just
as Jack polynomials specialize to classical bases of symmetric
function theory {\cite{DLM3}. Moreover, using a
$\beta$-generalization of the results of this article, a purely
combinatorial definition for these Jack polynomials in superspace
can remarkably be obtained (see \cite{DLM8b} for these
developments).

     The article is organized as follows. Section 2 first introduces the concept of superpartition.
     Then relevant results concerning the Grassmann algebra and
   symmetric superpolynomials are reviewed.  A simple interpretation
     of the latter, in terms of differentials forms, is also given.
     This section also includes  the definition of
     monomials in superspace and a formula for their products.

     Our main results are presented in section 3. It contains the superspace analog of the classical
     elementary symmetric functions, completely symmetric functions  and
     power-sum bases. The generating function for each of them
is displayed.  We point out  at this stage an interesting connection
between superpolynomials and de Rham complexes of symmetric
$p$-forms. Determinantal
     expressions that generalize classical formulas relating the
basis elements are presented. Furthermore, orthogonality and duality
relations are also established.

As already indicated, this work concerns, to a large extent, a
generalization of symmetric function theory. In laying down its
foundation, we generalize a vast number of basic results from this
theory which can be found for instance in \cite{Mac} and
\cite{StanBook} (Chap. 7).  Clearly, the core of most of our
derivations is bound to be a variation around the proofs of these
older results. We have chosen not to refer everywhere to the
relevant ``zero-fermionic degree'' version of the stated results to
avoid overquoting. But we acknowledge our debt in that regard  to
these two classic references.

%====================================

\section{Foundations}

\subsection{Superpartitions}

We recall that a partition
$\lambda=(\lambda_1,\lambda_2,\ldots,\lambda_\ell)$ of $n$, also
written as $\la\vdash n$, is an ordered set of integers such that:
$\lambda_1\geq\lambda_{2}\geq\ldots\geq\la_{\ell}\geq 0$ and
$\sum_{i=1}^{\ell}\la_i=n$. A particular juxtaposition of two
partitions gives a superpartition.

\begin{definition} \cite{DLM1} A
{\it superpartition} $\Lambda$ in the $ m $-fermion sector is a
sequence of non-negative integers separated by a semicolon such that
the sequence before the semicolon is a partition with $m$ distinct
parts, and such that the remaining sequence is a usual partition.
That is,  \begin{equation}\label{defsuperpart} \Lambda
=(\Lambda_1,\ldots,\Lambda_m;\Lambda_{m+1},\ldots,\Lambda_{N})\,
,\end{equation} where $ \Lambda_i>\Lambda_{i+1}\geq 0$ for $ i=1,
\ldots m-1$ and $\Lambda_j \ge \Lambda_{j+1}\geq 0 $ for
$j=m+1,\dots,N-1$.
\end{definition}

Given $\Lambda=(\Lambda^a ;\Lambda^s)$,  the partitions $\Lambda^a$
and $\Lambda^s$ are respectively called the {\it antisymmetric and
the symmetric components} of $\Lambda$. (From now on,  superscripts
$a$ and $s$ refer respectively to strictly and weakly decreasing
sequences of non-negative integers.) The {\it bosonic degree} (or
simply degree) of $\La$ is $|\Lambda|=\sum_{i=1}^{N}\Lambda_i$,
while its {\it fermionic degree} (or sector) is
$\overline{\underline{\Lambda}} = m$. Note that, in the zero-fermion
sector, the semicolon is usually omitted and $ \Lambda $  reduces
then to $\La^s$.

We say that the ordered set $\Lambda$ in (\ref{defsuperpart}) is a
superpartition of $(n|m)$ (a superparition of degree $n$ in the
fermionic sector $m$) if $|\Lambda|=n$ and
$\overline{\underline{\Lambda}} = m$, and write $\La\vdash(n|m)$.
The set composed of all superpartitions of $(n|m)$ is denoted
$\mbox{SPar}(n|m)$. When the fermionic degree is zero, we recover
standard partitions: $\mbox{SPar}(n|0)=\mbox{Par}(n)$.

We also define
\begin{equation} \mbox{SPar}(n):=\bigcup_{m\geq
0}\mbox{SPar}(n|m)\quad\mbox{and}\quad \mbox{SPar}:=\bigcup_{m,n\geq
0}\mbox{SPar}(n|m)\, ,\end{equation} with $
\mbox{SPar}(0|0)=\emptyset$ and $\mbox{SPar}(0|1)=\{(0;0)\}$. For
example, we have
\begin{equation}\mbox{SPar}(3|2)=\{\,(3,0;0),\,
(2,1;0),\,(2,0;1)\,(1,0;2),\,(1,0;1,1)\, \}\, .\end{equation} Notice
that $\mathrm{SPar}(n|m)$ is empty for all $n<m(m-1)/2$.

The length of a superpartition is
\begin{equation}{\ell}(\La):=\overline{\underline{\Lambda}}+{\ell}(\Lambda^s)
\quad\mbox{with}\quad {\ell}(\Lambda^s):=\mbox{Card}\{\La_i\in\La^s
\,:\,\La_i>0\}\, .\end{equation}
     With this
definition, ${\ell}\bigl((1,0;1,1)\bigr)=2+2=4$ (i.e., a zero-entry
in $\La^a$ contributes to the length of $\La$). To every
superpartition $\La$, we can also associate a unique partition
$\La^*$ obtained by deleting the semicolon and  reordering its parts
in non-increasing order.  For instance,
\begin{equation}
(5,2,1,0;6,5,5,2,2,1)^*=(6,5,5,5,2,2,2,1,1,0)=(6,5,5,5,2,2,2,1,1)\,.
\end{equation}
From this, we can introduce another notation for superpartitions. A
superpartition $\La=(\La^a;\La^s)$ can be viewed as the partition
$\Lambda^{*}$ in which every part of $\La^a$ is circled.
 If an entry $b$ of ${\La^a}$ also occurs in
$\La^s$, then we circle the leftmost $b$ appearing in $\La^*$. We
shall use $C[\Lambda]$ when refering to such a circled partition.
   For instance, \begin{equation}\Lambda=(3,1,0;4,3,2,1) \iff
C[\La]=(4,\cercle{3},3,2,\cercle{1},1,\cercle{0}).\end{equation}

\vskip0.2cm

The notation $C[\La]$ allows us to introduce a diagrammatic
representation of superpartitions. To each $\La$, we associate a
diagram, denoted by $D[\La]$.  It is obtained by first drawing the
Ferrer's diagram associated to $C[\La]$, that is, by drawing  a
diagram with ${C[\La]}_1$ boxes in the first row, ${C[\Lambda]}_2$
boxes in the second row and so forth, all rows being left justified.
In addition, if the $j$-th entry of ${C[\La]}$  is circled, then we
add a circle at the end of the $j$-th row of the diagram. We shall
further denote by ${\rm sh}(D[\Lambda])$ the {\it shape} of $D[\La]$
(including the circles). For example, with
$\Lambda=(3,1,0;4,3,2,1)$, we have as mentionned
$C[\La]=(4,\cercle{3},3,2,\cercle{1},1,\cercle{0})$, and thus
\begin{equation}
    D[\Lambda]=
{\tableau[scY]{&&&\\&&&\bl\tcercle{}\\&&\\&\\&\bl\tcercle{}\\&\bl\\\tcercle{}\bl\\
}} \, ,\end{equation} giving that ${\rm
sh}(D[\Lambda])=(4,4,3,2,2,1,1)$.

 The {\it conjugate} of a superpartition $\La$,
denoted by $\La'$, is the superpartition whose diagram is the
transposed (with respect to the main diagonal) of that of
$D[\Lambda]$. Hence, $(3,1,0;4,3,2,1)'=(6,4,1;3)$ since the
transposed of the previous diagram is

\begin{equation}
{\tableau[scY]{&&&&&&\bl\tcercle{}\\&&&&\bl\tcercle{}\\&&\\&\bl\tcercle{}\\
}}\end{equation}

%Recall that, for any partition $\la$, the
%conjugation can be defined by ${\la'}_j=\mbox{Card}\{k\,:\,\la_k\geq
%j\}$. Thus, in symbols, the conjugation of $\La\in\mbox{SPar}(n|m)$
%reads \begin{equation}\La'=({\La'}^a ; {\La'}^s)\, ,\end{equation}
%where\begin{equation}{\La'}^s=\left({\La^*}'\setminus{\La'}^a\right)^+\aand
%{{\La'}^a}_j=\mbox{Card}\{k\, :
%\,\La_k>\La_{m+1-j}\}\,,\end{equation} with $\la^+$ standing  for the
%partition obtained by reordering the parts of $\la$ non-increasingly.
\noindent Obviously, the conjugation of any superpartition $\La$
satisfies
\begin{equation}(\La')'=\La\quad\mbox{and}\quad(\La^*)'=(\La')^*\,
.\end{equation}

\begin{remark}
The description of superpartitions using  Ferrer's diagram with some
rows  ending with a circle makes clear that superpartitions are
equivalent to standard MacMahon diagrams, which are Ferrer's
diagrams with some corner cells marked (see for instance \cite{pak},
section 2.1.3). We shall nevertheless keep refering to
superpartition as superpartitions to be consistent with our previous
articles.

Furthermore, the notation $C[\Lambda]$ for a superpartition gives
immediately that the overpartitions introduced recently in
\cite{CoLo}  are special cases of superpartitions. Indeed,
overpartitions are circled superpartitions (with the circle replaced
by an overbar) that do not contain a possible circled zero. If we
denote by $s_N(n|m)$ the number of superpartitions of $(n|m)$ such
that ${\ell}(\La)\leq N$, then this connection implies that their
generating function is
\begin{equation}\label{gefct}
 \sum_{n,m, p\geq 0} s_{m+p} (n|m)
\,z^m y^{p} q^n=\frac{(-z;q)_\infty}{ (yq;q)_\infty} \quad {\rm with
}\quad (a;q)_\infty:=\prod_{n\geq 0}(1-aq^n)
\end{equation}
\end{remark}

To complete this subsection, we consider the natural ordering on
superpartitions.  We will first define it in terms of the Bruhat
order on compositions, and then later, in Corollary~\ref{coropart},
give a simpler characterization. Recall that a composition of $n$ is
simply a sequence of non-negative integers whose sum is equal to
$n$; in symbols $\mu=(\mu_1,\mu_2,\ldots)\in \mathrm{Comp}(n)$ iff
$\sum_i\mu_i =n$ and $\mu_i\geq0$ for all $i$.  The Bruhat ordering
on compositions is defined as follows. Given a composition $\la$, we
let $\la^+$ denote the partition obtained by reordering its parts in
non-increasing order. Now, $\la$ can be obtained from $\la^+$ by a
sequence of permutations. Among all permutations $w$ such that
$\lambda= w \lambda^+$, there exists a unique one, denoted
$w_{\lambda}$, of minimal length. For two compositions $\la$ and
$\mu$, we say that $\lambda \geq \mu $ if either $\lambda^+ > \mu^+$
in the usual dominance ordering or $\lambda^+=\mu^+$ and
$w_{\lambda} \leq w_{\mu}$ in the sense that the word $w_{\lambda}$
is a subword of $w_{\mu}$ (this is the Bruhat ordering on
permutations of the symmetric group). Recall that for two partitions
$\lambda$ and $\mu$ of the same degree, the dominance ordering is:
$\lambda \geq \mu$ iff
${\lambda}_1+\ldots+{\lambda}_k\geq{\mu}_1+\ldots+{\mu}_k$ for all
$k$.

Let $\La$ be a superpartition of $(n|m)$.  Then, to $\La$ is
associated  a unique composition of $n$, denoted by $\La^c$,
obtained by replacing the semicolon in $\Lambda$ by a comma.  We
thus have $\mathrm{Spar}(n)\subset\mathrm{Comp}(n)$, which leads to
a natural Bruhat ordering on superpartitions.

\begin{definition} \cite{DLM3} \label{bruhatorder}Let $\La,\,\Om\,\in\,\spar$.
The {\it Bruhat order}, denoted by $\leq$, is defined such that
$\Omega\leq\Lambda$ iff $\Om^c\leq\La^c$.
\end{definition}

For later purposes, we shall divide the Bruhat order into two
orders, depending on whether or not the superpartitions reorder to
the same partitions.
\begin{definition}\label{storder}Let $\La,\,\Om\,\in\,\spar$.
The {\it $S$ and $T$ orders} are respectively defined as follows:
\begin{equation}\begin{array}{lllllll}
\Omega\leq_S\Lambda&\mbox{if either}&\Om=\La&\mbox{or}&
\Om^*<\La^*\, ,&\\
\Omega\leq_T\Lambda&\mbox{if either}&\Om=\La&\mbox{or}& \Om^*=\La^*
&\mbox{and}&\Om^c<\La^c\, .\end{array}\end{equation}
\end{definition}
\n In order to describe other characterizations of these orders, we
need the following operators on compositions (or superpartitions):
\begin{equation}\label{defopST}\begin{array}{ll}
S_{ij}(\ldots,\la_i,\ldots,\la_j,\ldots)=&\ds\left\{\begin{array}{ll}
(\ldots,\la_i-1,\ldots,\la_j+1\ldots)&\mbox{if}\quad\la_i-\la_j>1\,
,\\
(\ldots,\la_i,\ldots,\la_j,\dots) &\mbox{otherwise}\,
,\end{array}\right.\\&\\
T_{ij}(\ldots,\la_i,\ldots,\la_j,\ldots)=&\ds\left\{\begin{array}{ll}
(\ldots,\la_j,\ldots,\la_i,\ldots)&\mbox{if}\quad\la_i-\la_j>0\,
,\\
(\ldots,\la_i,\ldots,\la_j,\dots) &\mbox{otherwise}\,.
\end{array}\right.
\end{array}\end{equation}

\begin{remark}
The $S$ order is precisely the ordering introduced in \cite{DLM1}.
It differs from the more precise ordering of \cite{DLM2}, which was
called $\leq^s$. In \cite{DLM3}, it is called the $h$ ordering.  See
also appendix B of \cite{DLM4}.
\end{remark}

The next lemma will lead to a simpler characterization of the order
on superpartitions.

\begin{lemma}\label{lemmaST}\cite{Mac, manivel}  Let $\la$ and $\om$
be two compositions of $n$. Then, $\la^+>\om^+$ iff there exists a
sequence $\{S_{i_1,j_1}, \ldots, S_{i_k,j_k}\}$ such that
\begin{equation} \om^+=S_{i_1,j_1}\cdots S_{i_k,j_k}\la^+\, .\end{equation}
Similarly, $\la^+=\om^+$ and $\la>\om$ iff there exists a sequence
$\{T_{i_1,j_1}, \ldots, T_{i_k,j_k}\}$ such that
\begin{equation} \om=T_{i_1,j_1}\cdots T_{i_k,j_k}\la\, .\end{equation}
\end{lemma}

Since a $T_{i,j}$ operation on a superpartition $\Lambda$ amounts to
removing the circle from row $i$ in $D[\Lambda]$ and adding it to
row $j$, the second part of this lemma can be translated for
superpartitions as: $\La>_T\Om$ iff ${\rm sh}(D[\La]) > {\rm
sh}(D[\Om])$ in the dominance order, where we recall that ${\rm
sh}(D[\La])$ is the shape of $D[\La]$ with the circles included.
This provides a simpler way to understand the order on
superpartitions.

\begin{corollary}  \label{coropart}
Let $\La,\,\Om\,\in\,\spar$.  We have that $\Omega\leq\Lambda$ iff
$\Omega^* < \Lambda^*$ or $\Omega^* = \Lambda^*$ and ${\rm
sh}(D[\Omega]) \leq {\rm sh}(D[\Lambda])$.
\end{corollary}

At this stage, we are in a position to establish the fundamental
property relating conjugation and Bruhat order which is that
  the Bruhat order is anti-conjugate (in the sense of the following proposition).

\begin{proposition}\label{uty}
Let $\La,\,\Om\, \in\, \spar (n|m)$. Then \begin{equation}\La\geq\Om
\quad\Longleftrightarrow\quad \Om'\geq\La'\,
.\end{equation}\end{proposition}
\begin{proof}It suffices to prove the result for the $S$ and $T$
orderings. The case $\La>_S\Om$, that is, $\La^*>\Om^*$, is a
well-known result on partitions (see for instance (1.11) of
\cite{Mac}).
 From Corollary~\ref{coropart}, the case $\La>_T\Om$
follows from the same argument.
\end{proof}

\vskip0.2cm
\begin{remark}Notice that we had before introduced
 as an alternative ordering
the dominance ordering on superpartitions,  denoted by $\leq_D $ and
defined as follows: $\Omega\leq_D\Lambda$ if either
$\Omega^*<\Lambda^*$ or $\Omega^*=\Lambda^*$ and
${\Omega}_1+\ldots+{\Omega}_k\leq {\Lambda}_1+\ldots{\Lambda}_k,$
$\forall\,\, k\, .$
 The
  usefulness of this ordering in special contexts lies in its simple
description in terms of inequalities. However, it is not the proper
generalization of the dominance order on partitions. In fact, it is
not as strict as  the  Bruhat ordering (i.e., more superpartitions
are comparable in this order than in the Bruhat ordering). This
follows from the second property of Lemma \ref{lemmaST} which
obviously implies that for superpartitions, the Bruhat ordering  is
a weak subposet of the dominance ordering, that is, $ \La
\geq\Om\,\Rightarrow\,\La \geq_D\Om$. However the converse is not
true.   For instance, if $\La=(5,2,1;4,3,3)$ and
$\Om=(4,3,0;5,3,2,1)$ we easily verify that $\La>_D \Om$. But since
$\Lambda^*=\Om^*$, ${\rm sh}(D[\La])=(6,4,3,3,3,2)$ and ${\rm
sh}(D[\Om])=(5,5,4,3,2,1,1)$ we have that $\La \not > \Om$ by
Corollary~\ref{coropart}. This corrects a loose implicit statement
in \cite{DLM3} concerning the expected equivalence of these two
orderings.
\end{remark}

\subsection{Ring of symmetric polynomials in superspace}

Let $\mathscr{B}=\{B_j\}$ and $\mathscr{F}=\{F_j\}$ be the formal
and  infinite sets composed of all bosonic (commutative) and
fermionic (anticommutative) quantities respectively.
% \begin{equation}[B_j,B_k]=[B_j,F_k]=\{F_j,F_k\}=0\, ,\end{equation}
% where $[A,B]$ and
% $\{A,B \}$ stand respectively for $AB-BA$ and $AB+BA$.
Thus, $\mathscr{S}=\mathscr{B}\oplus\mathscr{F}$ is
$\mathbb{Z}_2$-graded over any ring $\mathbb{A}$ when we identify
$\,^0\!\mathscr{S}$ with $\mathscr{B}$ and $\,^1\!\mathscr{S}$ with
$\mathscr{F}$.
% The degree
% of any element $s$ of $\mathscr{S}$, written $\hat{\pi}(s)$, is
% defined via \begin{equation}\hat{\pi}(s)=\left\{\begin{array}{ll}
% 0,& s\in \mathscr{B},\cr 1, & s\in \mathscr{F}.
% \end{array}\right.\end{equation}
% Consequently,
$\mathscr{S}$ possesses a linear endomorphism $\hat{\Pi}$, called
the parity operator,
 defined by
\begin{equation}\hat{\Pi}(s)=(-1)^{\hat{\pi}(s)} \, , \qquad {\rm where} \qquad \hat{\pi}(s)=\left\{\begin{array}{ll}
0,& s\in \mathscr{B},\cr 1, & s\in \mathscr{F}.
\end{array}\right.\end{equation}
In other words, the product of two bosons gives a boson, the product
of a boson and a fermion gives a fermion, and the product of two
fermions gives a boson.

An example of such a structure, is the Grassmann algebra over a
unital ring $\mathbb{A}$, denoted $\mathscr{G}_M(\mathbb{A})$.  It
is the algebra with identity $1\in\mathbb{A}$ generated by the $M$
anticommuting elements $\theta_1,\ldots,\theta_M $. We shall need
the following linear involution on the Grassmann algebra defined by
:\begin{equation}\overleftarrow{\,\theta_{j_1}\cdots\theta_{j_m}\,}:=\overrightarrow{\,\theta_{j_m}\cdots\theta_{j_1}\,}\,
\qquad {\rm where} \qquad
\overrightarrow{\,\theta_{j_1}\cdots\theta_{j_m}\,}:=\theta_{j_1}\cdots\theta_{j_m}\,
. \end{equation}
    In words, the operator $\overleftarrow{\phantom{abc}}$ reverses the order
of the anticommutative variables while
$\overrightarrow{\phantom{abc}}$  is simply the identity map. (The
explicit use of $\overrightarrow{\phantom{abc}}$ is not essential,
but it will make many formulas more symmetric and transparent.)
Using
 induction, we  get
\begin{equation}\overleftarrow{\,\theta_{j_1}\cdots\theta_{j_m}\,}=(-1)^{m(m-1)/2}\,
\overrightarrow{\,\theta_{j_1}\cdots\theta_{j_m}\,}\,.\end{equation}
    This result immediately  implies the following simple properties.

\begin{lemma}\label{ordredestheta} Let $\{\theta_1,\ldots,\theta_N\}$
and $\{\phi_1,\ldots,\phi_N\}$ be two sets of Grassmannian
variables.  Then \begin{equation} (\theta_{j_1}\phi_{j_1})\cdots
(\theta_{j_m}\phi_{j_m})=\overleftarrow{(\theta_{j_1}\cdots\theta_{j_m})}\overrightarrow{(\phi_{j_1}\cdots\phi_{j_m})}
=\overrightarrow{\,(\theta_{j_1}\cdots\theta_{j_m})\,}\overleftarrow{\,(\phi_{j_1}\cdots\phi_{j_m})\,}\,
\end{equation} and
    \begin{equation}
\overleftarrow{\overleftarrow{(\theta_{j_1}\cdots\theta_{j_m})}\overrightarrow{(\phi_{j_1}\cdots\phi_{j_m})}}
=\overrightarrow{(\phi_{j_1}\cdots\phi_{j_m})}\overleftarrow{(\theta_{j_1}\cdots\theta_{j_m})}\,
.\end{equation}\end{lemma}

Now, let $x=\{x_1,\ldots,x_N\}\subset\mathscr{B}$ and
$\theta=\{\theta_1,\ldots,\theta_N\}\subset\mathscr{F}$. We shall
let $\mathscr{P}( \mathbb{A})$ be the Grassmann algebra
$\mathscr{G}_N$ over the ring of polynomials in $x$ with
coefficients in $\mathbb A$. Note that $\mathscr{P}( \mathbb{A})$
can simply be considered as the ring of polynomials in  the
variables $x$ and $\theta$ over $\mathbb A$.

It is obvious that $\mathscr{P}$ is bi-graded with respect to the
bosonic and fermionic degrees, that is,
\begin{equation}\mathscr{P} =\bigoplus_{n,m\geq0}\mathscr{P}_{(n|m)}
\, ,\end{equation} where $\mathscr{P}_{(n|m)}$ is the finite
dimensional module made out of all homogeneous polynomials
$f(x,\theta)$ with degrees $n$ and $m$ in $x$ and $\theta$,
respectively.  Every  polynomial $f(x,\theta)$ in $\mathscr{P}$ also
possesses a bosonic and a fermionic part, i.e.,
% \begin{equation}\label{polxtheta}f(x,\theta)= \,^0\!f(x,\theta)+
% \,^1\!f(x,\theta)\, ,\end{equation}
$f(x,\theta)= \,^0\!f(x,\theta)+
 \,^1\!f(x,\theta) $ where $\,^0\!f(x,\theta) \in
\mathscr{B}$ and $\,^1\!f(x,\theta) \in \mathscr{F}$. We have that
$\,^0\!f(x,\theta)$ consists of the monomials of $f$ with an even
degree in $\theta$ while $\,^1\!f(x,\theta)$ consists of those
monomials with an odd degree in $\theta$. Purely fermionic
polynomials (i.e., elements of $\mathscr{P} _{(n|m)}$ with $m$ odd)
have some nice properties. As an example, consider the following
proposition that shall be useful in the subsequent sections.

\begin{proposition}\label{expfermions}Let
$\tilde{f}=\{\tilde{f}_0,\tilde{f}_1,\ldots\}$ and
$\tilde{g}=\{\tilde{g}_0,\tilde{g}_1,\ldots\}$ be two sequences of
fermionic polynomials parametrized by non-negative integers.  Let
also
\begin{equation}\tilde{f}_{\mu}:=\tilde{f}_{\mu_1}\tilde{f}_{\mu_2}\cdots\aand
\tilde{g}_{\mu}:=\tilde{g}_{\mu_1}\tilde{g}_{\mu_2}\cdots
\end{equation}where $ \mu$ belongs to $
\mathrm{Par}_a(n)$, the set of partitions of $n$ with strictly
decreasing parts.    Then
\begin{equation}\exp\left[\sum_{n=0
}^{N-1}\tilde{f}_n\,\tilde{g}_n\right]=\sum_{n=0}^{N(N-1)/2}\sum_{\mu\in\mathrm{Par}_a(n)}
\overleftarrow{\phantom{\Big|}\tilde{f}_\mu\phantom{\Big|}}\,\overrightarrow{\phantom{\Big|}\tilde{g}_\mu\phantom{\Big|}}\,
.\end{equation}
\end{proposition}
\begin{proof} Due to the fermionic character of $\tilde{f}$ and
$\tilde{g}$ ($\tilde{f}^2=\tilde{g}^2=0$ for instance), we have
\begin{equation}\begin{array}{lcl}\exp\left[\sum_{n=0
}^{N-1}\tilde{f}_n\,\tilde{g}_n\right]&=&\ds \prod_{0\leq  n\leq
N-1} (1+\tilde{f}_n\,\tilde{g}_n)\\&=&\ds 1+\sum_{0\leq n\leq
N-1}\tilde{f}_n\,\tilde{g}_n+\sum_{0\leq m<n\leq
N-1}\tilde{f}_m\,\tilde{g}_m
\tilde{f}_n\,\tilde{g}_n+\ldots\end{array}\end{equation} Since every
term in the last equality can be reordered by Lemma
\ref{ordredestheta}, the proof follows.\end{proof}

    We finally define what we consider as
symmetric polynomials in superspace.   The algebra of symmetric
superpolynomials over the ring $\mathbb{A}$, denoted by
$\mathscr{P}^{S_N}(\mathbb{A})$ or by $\mathbb{A}[x_1,\ldots,x_N,
\theta_1,\ldots,\theta_N]^{S_N}$, is a subalgebra of $\mathscr{P}$.
As mentioned in the introduction, $\mathscr{P}^{S_N}$ is made out of
all $f(x,\theta)\in\mathscr{P}$ invariant under the diagonal action
of the symmetric group $S_N$ on the two sets of variables.

To be more explicit, we introduce $K_{ij}$ and $\kappa_{ij}$, two
distinct polynomial realizations of the transposition $(i,j)\in
S_N$:\begin{equation}K_{ij}f(x_i,x_j,\theta_i,\theta_j)=f(x_j,x_i,\theta_i,\theta_j)
\, ,\qquad
\kappa_{ij}f(x_i,x_j,\theta_i,\theta_j)=f(x_i,x_j,\theta_j,\theta_i)\,
,\end{equation} for all $f\in \mathscr{P}$.  Since every permutation
is generated by products of elementary transpositions $(i,i+1)\in
S_N$, we can define symmetric superpolynomials as follows.
\begin{definition} \label{defsupersym}
A  polynomial $f(x,\theta)\in\mathscr{P}$ is {\it symmetric} if and
only if \begin{equation}\mc{K}_{i,i+1}
f(x,\theta)=f(x,\theta)\quad\mbox{where}\quad\mc{K}_{i,i+1}:=\kappa_{i,i+1}K_{i,i+1}\end{equation}
for all $i\in\{1,2,\ldots,N-1\}$.\end{definition}

    \n Since every
monomial $\theta_J=\theta_{j_1}\cdots\theta_{j_m}$ is completely
antisymmetric,
% , that is, \begin{equation}
% \kappa_{ik}\theta_J=\left\{\begin{array}{ll}
% -\,\theta_J,&\mbox{if}\quad i\,, k\,\in J,\\
% \phantom{-}\,\theta_J,&\mbox{if}\quad i\,, k\,\not\in
% J.\end{array}\right.\end{equation}
  %   This
we have the following result.

\begin{lemma}\label{antisymsym}Let $f(x,\theta)\in \mathscr{P}$ be
expressed as:
\begin{equation} f(x,\theta)=\sum_{m\geq 0}\,\sum_{{1\leq j_1<\ldots<j_m\leq
N}}\,f^{j_1,\ldots,j_m}(x)\,\theta_{j_1}\cdots\theta_{j_m}\, .
\end{equation} If  $f(x,\theta)$ is symmetric, then each polynomial
$f^{j_1,\ldots,j_m}(x)$ is completely antisymmetric in the set of
variables $y:=\{x_{j_1},\ldots,x_{j_m}\}$ and completely symmetric
in the set of variables $x\setminus y$.\end{lemma}

 \begin{remark} The symmetric superpolynomials are completely different
from the ``supersymmetric polynomials'' previously considered in the
literature. Recall that what is  called a supersymmetric polynomial
(see e.g., \cite{Stem}) is first of all a doubly symmetric
polynomial in two distinct sets of ordinary (commuting) variables
$x_1,\cdots x_m$ and $y_1,\cdots, y_n$, i.e.,  invariant under
independent permutations of the $x_i$'s and the $y_i$'s. It is said
to be supersymmetric if, in addition, it satisfies the following
cancelation condition: by substituting $x_1=t$ and $y_1=t$, the
polynomial becomes independent of $t$. An example of a generating
function for such polynomials is
\begin{equation}
\prod_{i=1}^m (1-qx_i)\prod_{j=1}^n(1-qy_j)^{-1}= \sum_{r\geq 0}
p_{(r)}(x,y) q^r
\end{equation}
This generating function is known to appear in the context of
classical Lie superalgebras (as a superdeterminant) \cite{Kac}.
Actually most of the work on supersymmetric polynomials is motivated
by its connection with superalgebras. For an example of such an
early work, see \cite{jarvis}. More precisions and references are
also available in \cite{brenti, moens}. The key differences between
these supersymmetric polynomials and our symmetric superpolynomials
should
     be clear. In our case, we symmetrize two sets of variables with
respect to the diagonal action of the symmetric group, with one of
the two sets being made out of Grassmannian variables.
\end{remark}

\subsection{Geometric interpretation of polynomials in superspace}

Symmetric functions can be interpreted as symmetric 0-forms $f$
acting on a manifold: $K_{ij}f(x)=f(x)$ where $x$ is a local
coordinate system.  Similarly, symmetric superfunctions in the
$p$-fermion sector can be reinterpreted as symmetric $p$-forms $f^p$
acting on the same manifold: ${\mathcal K}_{ij}f^p(x)=f^p(x)$.
Thus, the set of all symmetric superfunctions is in correspondence
with the completely symmetric de Rham complex.  This geometric point
of view is   briefly explained in this subsection. (Note that none
of our results relies on this observation.)

We consider a Riemannian manifold $\M$ of dimension $N$ with metric
$g_{ij}$ and  let $x=\{x^1,\ldots,x^N\}$ denote a coordinate system
on a given subset of $\M$.
  On the
tangent bundle, we choose an orthonormal coordinate frame $
\{\partial_{1},\ldots,\d_{N}\}$. As usual, $ \{dx^1,\ldots,dx^N\}$
denotes the dual basis that belongs to the cotangent bundle, i.e.,
$dx^i(\d_{x^{j}})=\delta^{i}_j$. The set of all $p$-form fields on
$\M$ is a vector space denoted by ${\bigwedge}^p$.  Each $p$-form
can be written as
\begin{equation}\label{developform}\alpha^p(x)=\sum_{1\leq j_1<\ldots<j_p\leq N}
\alpha_{j_1,\ldots,j_p}(x) dx^{j_1}\wedge\cdots\wedge dx^{j_p}\,  ,
\end{equation} where the exterior  (wedge)  product is antisymmetric
: $dx^i\wedge dx^j=-dx^j\wedge dx^i$. Let $d$ be the exterior
differentiation on forms, whose action is
\begin{equation}d\, \alpha^p(x)=\sum_{1\leq k,j_1,\ldots,j_p\leq
N}[\d_{x^k}\alpha_{j_1,\ldots,j_p}(x)]dx^k\wedge
dx^{j_1}\wedge\cdots\wedge dx^{j_p}\,  . \end{equation}
    This operation is used to
define the de Rham complex of $\M$:
\begin{equation}0\longrightarrow\mathbb{R}
\longrightarrow{\bigwedge}^0\stackrel{d}{\longrightarrow}
{\bigwedge}^1\stackrel{d}{\longrightarrow}\cdots\stackrel{d}{\longrightarrow}{\bigwedge}^N
\stackrel{d}{\longrightarrow}0\, .\end{equation}

In order to represent our Grassmannian variables $\theta^j$ and
$\theta^j \,^\dagger$  in terms of forms, we introduce the two
operators $\hat{e}_{dx^j}$ and $\hat{\i}_{\d_{x_k}}$,
   where $\hat{e}_\alpha$ and $\hat{\i}_{v}$ respectively stand for the
  left exterior product by the form $\alpha$ and
  the interior product (contraction) with respect to the vector field $v$.
  These operators satisfy a
  Clifford (fermionic) algebra
\begin{equation}\{\hat{e}_{dx^i}\,,\,\hat{\i}_{\d_{x_k}}\,\}=\delta^i_k\aand
\{\hat{e}_{dx^i}\,,\,\hat{e}_{dx^j}\}=0=\{\hat{\i}_{\d_{x_j}}\,,\,
\hat{\i}_{\d_{x_k}}\,\}\, .\end{equation} This implies that the
$\theta^j$'s and  $\theta^j\,^\dagger$, as operators,
    can be realized as follows:
  \begin{equation} \theta^j\sim\hat{e}_{dx^j}\aand \theta^j
  \,^\dagger \sim g^{jk}\,\hat{\i}_{\d_{x_k}}
   \,,\end{equation}
that is,
\begin{equation} dx^j\sim \theta^j(1)\aand\theta^j\,^\dagger\sim
g^{jk}\,\partial_{\theta^k}\wwhere
\partial_{\theta^k}:=\frac{\partial\,}{\partial \theta^k}\, .\end{equation}
Note that introducing the Grassmannian variables as operators is
needed to enforce the wedge product of the forms $dx_j$. Moreover,
if $\alpha^p$ is a generic $p$-form field and $\hat{\pi}\,
:\,{\bigwedge}^p\rightarrow {\bigwedge}^p$ is the operator defined
by
\begin{equation} \hat{\pi}_p:=
\theta^j\partial_{\theta^j}=\theta^jg_{jk}\theta^k\,^\dagger\quad\Longrightarrow\quad
\hat{\pi}_p\alpha^p=p\,\alpha^p\,,  \end{equation}then
$\hat{\Pi}_p:=(-1)^{\hat{\pi}_p}$ is involutive. (Manifestly,
$\hat{\Pi}_p=\hat{\Pi}$, the parity operator introduced previously.)
This involution is also an isometry in the Hilbert space scalar
product.  The operator $\hat{\Pi}_p$ induces a natural
$\mathbb{Z}_2$ grading in the de Rham complex.

The construction of  the symmetric de Rham complex
%, denoted $\mathrm{SRham}$,
 is obtained as follows. We make a change of
coordinates: $x\rightarrow f(x)$, where
$f=\{f^n\}:=\{f^1,\ldots,f^N\}$ is an $N$-tuple of symmetric and
independent functions of $x$.  For instance, $f^n$ could be an
elementary symmetric function $e_n$, a complete symmetric function
$h_n$, or a power sum $p_n$ (see Section 3). This implies a change
of basis in the cotangent bundle: $dx\rightarrow df(x)$. Explicitly,
\begin{equation} df^n=\sum_i (\partial_i f^n)(x) \,dx_i\sim
\tilde{f}^n=\sum_i (\partial_i f^n)(x) \,\theta_i\, .\end{equation}
In other words, $d{f}$ is a new set of  ``fermionic'' variables
invariant under any permutation of the $x_j$'s.
% \footnote{Of course,
% this change of basis is well defined in $U\subset\M$ if and only if
% the Jacobian determinant $J(f,x)$ is not zero.  For any standard
% basis of the symmetric function space (i.e., the $e_n$'s, $h_n$'s or
% $p_n$'s), we easily verify that the Jacobian is, up to a sign, equal
% to the Vandermonde determinant.   Thus, the
% change of coordinates is well defined for any locus $U$  such as the
% following one: $L_U:=\{x(p)\,:\,x^1(p)>x^2(p)>\ldots\, ;\,\forall\,
% p\in U\}\, .$}

  These remarks explicitly show that
   symmetric polynomials in superspace can be interpreted as
  symmetric differential forms.  We stress that the diagonal
   action of the symmetric group $S_N$ comes naturally in the
   geometric perspective.  Note finally
that for an Euclidian superspace (relevant to our context), the
position (upper or lower) of the indices does not matter.

\subsection{Monomial basis}

The monomial symmetric functions in superspace, denoted by
$m_{\Lambda}= m_{\Lambda}(x,\theta)$, are  the superanalog of the
monomial symmetric functions.
\begin{definition} \cite{DLM1} \label{defsupermono}To each
$\La\in\mathrm{SPar}(n|m)$, we associate the monomial symmetric
function
%$m_{\Lambda}\in\mathscr{P}^{S_N}_{(n|m)}$
%, called {\it symmetric supermonomial},
%defined by
\begin{equation}
m_{\Lambda}={\sum_{\sigma\in S_{N}}}' \theta^{\sigma(1, \ldots,
m)}x^{\sigma(\Lambda)} \, ,
\end{equation} where the prime indicates that the  summation is
restricted to distinct terms, and where
\begin{equation}
x^{\sigma(\Lambda)}=x_1^{\Lambda_{\sigma(1)}} \cdots
x_m^{\Lambda_{\sigma(m)}} x_{m+1}^{\Lambda_{\sigma(m+1)}} \cdots
x_{N}^{\Lambda_{\sigma(N)}}\quad\mbox{and}\quad\theta^{\sigma(1,
\ldots, m)} = \theta_{\sigma(1)} \cdots \theta_{\sigma(m)} \, .
\end{equation}\end{definition}

\noindent Obviously, the previous definition can be replaced by the
following:
\begin{equation} \label{monod}
         m_{\Lambda}= \frac{1}{n_{\Lambda}!} \sum_{\sigma \in S_N}
\mc{K}_{\sigma}
         \left( \theta_1 \cdots \theta_m x^{\Lambda}\right) \end{equation}
    with
\begin{equation} \label{augmented}
 n_{\Lambda}! =n_{\Lambda^s}!:=
n_{\Lambda^s}(0)!\, n_{\Lambda^s}(1)! \,
         n_{\Lambda^s}(2) ! \cdots \, ,
\end{equation}
    where $n_{\Lambda^s}(i)$ indicates  the number of $i$'s in
$\Lambda^s$,
         the symmetric part of $\Lambda=(\Lambda^a;\Lambda^s)$.
Moreover, $\mc{K}_{\sigma}$ stands for $\mc{K}_{i_1,i_1+1} \cdots
\mc{K}_{i_n,i_n+1}$ when the element $\sigma$ of the symmetric group
$S_N$ is written  in terms of elementary transpositions, i.e.,
$\sigma = \sigma_{i_1} \cdots \sigma_{i_n}$.  Notice that the
monomial symmetric function $m_\La$, with $\La\vdash(n|m)$,  belongs
to $\mathscr{P}_{(n|m)}(\mathbb{Z)}$, the space of superpolynomials
of degree $(n|m)$ with integer coefficients.

\begin{theorem}
The set $\{m_\La\}_{\La\vdash(n|m)}:=\{m_\La \,:\,
\La\in\spar(n|m)\}$ is a basis of
$\mathscr{P}^{S_N}_{(n|m)}(\mathbb{Z})$.
\end{theorem}
\begin{proof}
Each  polynomial $f(x,\theta)$ of degree
 $(n|m)$, with $N$ variables and with integer coefficients, can be
 expressed as a sum of monomials of the type
 $ \theta_{j_1} \dots \theta_{j_m} x^{\mu}$,
 with coefficient $a_{\mu}^{j_1,\dots,j_m} \in \mathbb Z$,
 and where $\mu$ is a composition of $N$.  Let $\Omega^a$ be the
 reordering of the entries
 $({\mu_{j_1},\dots,\mu_{j_m}})$, and let
 $\Om^s$ be the reordering of the remaining entries of $\mu$.
 Because the  polynomial $f(x,\theta)$ is also symmetric,
 $f(x,\theta)$
 is by definition invariant under the action of $\mathcal K_{\sigma}$,
 for any $\sigma \in S_N$.   Therefore,
 $a_{\mu}^{j_1,\dots,j_m}$
 must be equal, up to a sign, to the coefficient
 $a_{\Omega}^{1,2,\dots,m}$ of
 $\theta_1 \cdots \theta_m x^{\Om}$ in $f(x,\theta)$,
 where $\Om=(\Om^a;\Om^s)$.
 Note that  from Lemma~\ref{antisymsym}, $\Omega^a$ needs to have
 distinct parts,
 which means that $\Omega$ is a superpartition.
 This gives that $f(x,\theta)-a_{\Omega}^{1,2,\dots,m} m_{\Omega}$
 does not contain any monomial that is also a monomial of $m_{\Om}$,
 since otherwise
 it would need by symmetry to contain the monomial $\theta_1 \cdots
 \theta_m x^{\Om}$.

 Now, consider any total order on superpartitions, and let $\La$ be
 the highest superpartition in this order such that there is a
 monomial of $m_{\La}$ appearing in $f(x,\theta)$.
   By the previous argument, $f(x,\theta)- a_{\Lambda}^{1,\dots,m}m_{\La}$ is
 a symmetric superpolynomial such that no monomial belonging to
 $m_{\La}$ appears
 in its expansion.  Since no monomial of $m_{\La}$ appears in any other
 monomial of $m_{\Om}$, for $\Om\neq \La$, the proof follows by induction.
\end{proof}

\begin{corollary} The set $\{m_\La\}_{\La}:=\{m_\La \,:\,
\La\in\spar\}$ is a basis of $\mathscr{P}^{S_N}(\mathbb{Z})$.
\end{corollary}

This corollary implies that $\mathscr{P}^{S_N}(\mathbb{Z})$ could
also be defined  as the free $\mathbb{Z}$-module spanned by the set
of monomial symmetric functions in superspace.

To end this section, we give a formula for the expansion
coefficients of the product of two monomial symmetric functions
 in terms of monomial symmetric functions.
%This furnishes an illustrative example of  what
%could be called supercombinatorics.
In this kind of calculation, the standard counting of combinatorial
objects  is affected by signs resulting from the reordering of
fermionic variables.

\begin{definition} \label{deffilling} Let $\La \in \spar(n|m)$,
$\Omega \in \spar(n'|m')$ and $\Gamma \in \spar(n+n'|m+m')$. In each
box or circle of $D[\La]$ , we write a letter $a$. In its $i$-{th}
circle (the one corresponding to $\Lambda_i$), we add the label $i$
to the letter $a$. We do the same process for  $D[\Omega]$ replacing
$a$ by $b$. We then define $\mathcal T[ \La,\Om;\Gamma]$ to be the
set of distinct fillings of $D[\Gamma]$ with the letters of $D[\La]$
and $D[\Omega]$ obeying the following rules:
\begin{enumerate}
\item the circles of $D[\Gamma]$ can only be filled with labeled
letters (an $a_i$ or a $b_i$);
\item each row of the filling of $D[\La]$ is reproduced
in a single and distinct row of the filling of $D[\Gamma]$; in other
words, rows of $D[\La]$ cannot be split and two rows of $D[\La]$
cannot be put within a single row of $D[\Gamma]$;
% that is,
% for all $i$, there is one row of the filling of $D[\Gamma]$ where
% unlabeled $a$'s appear exactly $\Lambda_i$ times, with an additional
% $a_i$ if $i\leq m$.  Furthermore, if $\Lambda$ has $j$ entries equal
% to $\Lambda_i$, then there are $j$ such rows in the filling of
% $D[\Gamma]$;
% \item each row of the filling of $D[\Omega]$ is reproduced
% in a single and distinct row of the filling of $D[\Gamma]$;
\item rule 2 also holds when $D[\La]$ is replaced by $D[\Omega]$;
\item in each row, the unlabeled $a$'s appear to the left of the
unlabeled $b$'s.
\end{enumerate}
\end{definition}

 For instance, there are three possible fillings of
$(2,1,0;1^3)$ with $(1,0;1)$ and $(0;2,1^2)$:
\begin{equation}\label{exfilling}
{\tableau[scY]{b&b&\bl\tcercle{$a_2$}\\a&\bl\tcercle{$a_1$}\\a\\b\\b\\\bl\tcercle{$b_1$}}}
\qquad
{\tableau[scY]{b&b&\bl\tcercle{$a_2$}\\a&\bl\tcercle{$a_1$}\\b\\a\\b\\\bl\tcercle{$b_1$}}}
\qquad
{\tableau[scY]{b&b&\bl\tcercle{$a_2$}\\a&\bl\tcercle{$a_1$}\\b\\b\\a\\\bl\tcercle{$b_1$}}}\,
.
  \end{equation}\label{exone}
    There are also three possible fillings of
$ (3,1,0;1^2)$ with $(1,0;1)$ and $(0;2,1^2)$:
\begin{equation}\label{exfilling2}
{\tableau[scY]{a&b&b&\bl\tcercle{$a_1$}\\b&\bl\tcercle{$a_2$}\\a\\b\\\bl\tcercle{$b_1$}}}
\qquad
{\tableau[scY]{a&b&b&\bl\tcercle{$a_1$}\\b&\bl\tcercle{$a_2$}\\b\\a\\\bl\tcercle{$b_1$}}}
\qquad
{\tableau[scY]{a&b&b&\bl\tcercle{$a_1$}\\a&\bl\tcercle{$b_1$}\\b\\b\\\bl\tcercle{$a_2$}}}
\,.\end{equation}

\begin{definition}\label{defewightfilling} Let
$T \in  \mathcal T[ \La,\Om;\Gamma]$, with
$\overline{\underline{\Lambda}}=m$ and
$\overline{\underline{\Om}}=m'$. The weight of $T$, denoted by
$\hat{w}[T]$, corresponds to the sign of the permutation needed to
reorder the content of the circles in the filling of $D[\Gamma]$ so
that from top to bottom they read as $a_1 \dots a_m b_1 \dots
b_{m'}$.
\end{definition}

In the example (\ref{exfilling}), each term has weight $-1$ (odd
parity). The oddness of these fillings comes from the transposition
that is needed to reorder $a_1$ and $a_2$. In the second example
(\ref{exone}), the two first fillings are even while the last
filling is odd due to the needed transposition of $a_2$ and $b_1$.
As we shall see in the next proposition, the two previous sets lead
respectively to the coefficients of $m_{(2,1,0;1^3)}$ and
$m_{(3,1,0;1^2)}$ in the product of $m_{(1,0;1)}$ and
$m_{(0;2,1^2)}$, that is,
\begin{equation}
m_{(1,0;1)}\,m_{(0;2,1^2)}=\underbrace{(-3\,)}_{-1-1-1}\times\,m_{(2,1,0;1^3)}+\underbrace{(\,1\,)}_{+1+1-1}\times\,m_{(3,1,0;1^2)}+\mbox{other
terms}\, .\end{equation}
\begin{proposition}
Let $m_\La$ and $m_\Om$ be any two monomial symmetric functions in
superspace.  Then
\begin{equation}\label{propmono}m_\La \,m_\Om
=\sum_{\Gamma\in\mathrm{SPar}} N^\Gamma_{\La,\Om}\,m_\Gamma\,
,\end{equation} where the integer $N^\Gamma_{\La,\Om}=(-1)^{
\overline{\underline{\Lambda}}\, \cdot\, \overline{\underline{\Om}}}
\, N^\Gamma_{\Om,\La}$ is given by
\begin{equation}N^\Gamma_{\La,\Om}:= \sum_{T \in  \mathcal T[
\La,\Om;\Gamma] } \hat w[T]\, .\end{equation}
\end{proposition}
  \begin{proof}
 From the symmetry property in Definition~\ref{defsupermono}, the
coefficient $N^\Gamma_{\La,\Om}$ is simply given by the coefficient
of $\theta_{\{1,\ldots,m+p\}}x^\Gamma$ in $m_\Lambda m_{\Om}$.  The
terms contributing to this coefficient correspond to all distinct
permutations $\sigma$ and $w$ of the entries of $\Lambda$ and
$\Omega$ respectively such that
\begin{equation} \label{proofmono2}
\Gamma=
(\Lambda_{\sigma(1)}+\Omega_{w(1)},\dots,\Lambda_{\sigma(N)}+\Omega_{w(N)})
\, ,
\end{equation}
where the entries of $\Lambda^a$ and $\Omega^a$ are distributed
among the first $m+p$ entries (no two in the same position). But
this set is easily seen to be in correspondence with the fillings in
$\mathcal T[ \La,\Om;\Gamma]$ when realizing that labeled letters
simply give the positions of the fermions in $C[\Gamma]$ (the
circled version of $\Gamma$). The only remaining problem is thus the
ordering of the fermions.  In (\ref{proofmono2}), from the
definition of monomial symmetric functions, the sign of the
contribution is equal to the sign of the permutation needed to
reorder the fermionic entries of $\Lambda^a$ and $\Omega^a$ that are
distributed among the first $m+p$ entries so that they correspond to
$({\La^a}, \Om^a)$.  But this is simply the sign of the permutation
that reorders the circled entries in the corresponding filling of
$D[\Gamma]$ such that they read as $a_1 \dots a_m b_1 \dots b_p$.
\end{proof}

%====================
\section{Generating functions and multiplicative bases}

In the theory of symmetric functions, the number of variables is
usually irrelevant, and can be set for convenience to be equal to
infinity.  In a similar way, we shall consider from now on that,
unless otherwise specified, the number of $x$ and $\theta$ variables
is infinite, and denote the ring of symmetric superfunctions as
$\mathscr{P}^{S_\infty}$.

\subsection{Elementary symmetric functions}

Let $J=\{j_1,\ldots,j_r\}$ with $1\leq j_1<j_2 <j_3 \cdots$ and let
$\# J:=\mbox{Card} \, J$.  The  $n$-th {\it bosonic and fermionic
elementary symmetric  functions}, for $n \geq 1$,  are defined
respectively by
\begin{equation}e_n:=\sum_{J;\,\# J=n}x_{j_1}\cdots
x_{j_n}\qquad {\rm and}\qquad \tilde{e}_n:=\sum_{i\geq
1}\sum_{\substack{J;\,\# J=n\\ i\not\in J}}\theta_i\,x_{j_1}\cdots
x_{j_n}\end{equation}
%and
% \begin{equation}\tilde{e}_n:=\sum_{1\leq i\leq
% N}\sum_{\substack{J;\,\# J=n\\ i\not\in
% J}}\theta_i\,x_{j_1}\cdots x_{j_n}\in\mathscr{F}\, .\end{equation}
In addition, we impose \begin{equation}e_0=1\;\quad {\rm and } \quad
\tilde{e}_0=\sum_i\theta_i\, .
\end{equation}
     So, in terms of monomials, we have
\begin{equation}e_n=m_{(1^n)}\;, \qquad \tilde{e}_n=m_{(0;1^n)}\; .
\end{equation}

We introduce two parameters: $t\in\mathscr{B}$ and
$\tau\in\mathscr{F}$.  It is easy to verify that the generating
function for the elementary  symmetric functions is
\begin{equation}\label{defE} E(t,\tau):=\sum_{n=0}^\infty t^n(e_n
+\tau\tilde{e}_n) =\prod_{i=1}^\infty (1+t x_i + \tau \theta_i)\,
.\end{equation} Actually, to go from the usual generating function
$E(t):=E(t,0)$ to the new one, one simply replaces $x_i\rw
x_i+\tau'\ta_i$ and redefines $\tau'=\tau/t$, an operation that
makes manifest the invariance of $E(t,\tau)$ under the simultaneous
interchange of the $x_i$'s and the $\ta_i$'s.

   From an analytic point of view, the fermionic elementary
 symmetric functions are obtained by exterior differentiation:
\begin{equation}\tilde{e}_{n-1}(x,\theta)\sim\tilde{e}_{n-1}(x,dx)=d\,e_n(x)\,
,\end{equation}for all $n \geq 1$.  How can we explain that the
generating function (\ref{defE}) leads precisely to the  fermionic
elementary  functions that are obtained by the action of the
exterior derivative of the elementary symmetric function?  The
rationale for this feature turns out to be rather simple.  Indeed,
let $\tau:= \, t \, dt$ and define $D$ to act on a function $f(x,t)$
as a tensor-product  derivative:
\begin{equation} D f:= dt \wedge df \, .\end{equation}
In consequence, we formally have
\begin{equation}(1+t x_i +\tau \theta_i)\sim(1+D)(1+tx_i) \aand
E(t,\tau)\sim(1+D)\, E(t)\,,\end{equation} which is the desired
link.
%  The last relation explains why the fermionic
% elementary symmetric functions obtained from the generating function
% $E(t,\tau)$ are precisely represented by those resulting from  the
% action of $d$ on the bosonic ones.

In order to obtain a new basis of the symmetric superpolynomial
algebra, we associate, to each superpartition
$\La=(\La_1,\ldots,\La_m;\La_{m+1},\ldots, \La_\ell)$ of $(n|m)$, a
 polynomial $e_\La\in\mathscr{P}^{S_\infty}_{(n|m)}$ defined by
\begin{equation}e_\La:=\prod_{i=1}^m\tilde{e}_{\La_{i}}\prod_{j=m+1}^\ell
{e}_{\La_{j}}\, , \end{equation} Note that the product of
anticommutative quantities is always done from left to right:
$\prod_{i=1}^N F_i:=F_1F_2\cdots F_N$. We stress that the ordering
matters in the fermionic sector since for instance
\begin{equation}
{e}_{(3,0;4,1)}= {\tilde e}_{3}{\tilde e}_{0}{e}_{4}{e}_{1}= -
{\tilde e}_{0}{\tilde e}_{3}{e}_{4}{e}_{1}
 \end{equation}

\begin{theorem} \label{theoebase}Let $\La$ be a superpartition of
$(n|m)$ and $\La'$ its conjugate. Then \begin{equation}\label{einm}
\versg{ e_\La}=\tilde{e}_{\La_m}\cdots
\tilde{e}_{\La_1}{e}_{\La_{m+1}}\cdots{e}_{\La_{N}}=m_{\La'}+\sum_{\Om
< \La'} N_{\La}^{\Om}\,m_\Om\, ,\end{equation} where $N_{\La}^{\Om}$
is an integer. Hence, $\{\,e_\La \,:\, \La\vdash(n|m)\}$ is a basis
of $\mathscr{P}^{S_\infty}_{(n|m)}(\mathbb{Z})$.\end{theorem}
\begin{proof}
We first observe that $\versg{ e_\La}=(-1)^{m(m-1)/2} e_{C[\La]}$,
where $C[\La]$ denotes as usual the partition $\La^*$ in which
fermionic parts of $\La$ are identified by a circle. Then, assuming
that we work in $N$ variables, the monomials $\theta_J x^\nu$ that
appear in the expansion of $e_{C[\La]}$ are in correspondence with
the fillings of $D[\Lambda']$ with the letters $1,\dots,N$ such
that:
\begin{enumerate}
\item the non-circled entries in the filling of $D[\Lambda']$
increase when going
       down in a column;
\item if a column contains a circle, then the entry that fills the
circle cannot
       appear anywhere else in the column.
\end{enumerate}
The correspondence follows because the reading of the $i$-th column
corresponds to one monomial of $e_{\Lambda_i}$ (or $\tilde
e_{\Lambda_i}$). To be more specific, if the reading of the column
is $j_1,\dots,j_{\Lambda_i}$ (with a possible  extra letter $a$), it
corresponds to the monomial $x_{j_1} \cdots x_{j_{\Lambda_i}}$ (or
$\theta_a x_{j_1} \cdots x_{j_{\Lambda_i}}$). The first condition
ensures that we do not count the permutations of $x_{j_1} \cdots
x_{j_{\Lambda_i}}$ as distinct monomials.  The second one ensures
that in the fermionic case, the index of $\theta_a$ is distinct from
the index of the variables $x_{j_1},\dots,x_{j_{\Lambda_i}}$.

Now, to obtain the coefficient $N^{\Omega}_\Lambda$, it suffices to
compute the coefficient of $\theta_{C[\Omega]} x^{C[\Omega]}$ in
$e_{C[\Lambda]}$, where $\theta_{C[\Omega]}$ represents the product
of the fermionic entries of $C[\Omega]$ read from left to right.
Note that this coefficient has the same sign as
$\theta_{\{1,\dots,m\}}x^{\Omega}$ in $m_{\Om}$ and there is thus no
need to compensate by a sign factor. The monomials contributing to
$N^{\Omega}_\Lambda$ are therefore fillings of $D[\Lambda']$
(obeying the two conditions given above) with the letter $i$
appearing $C[\Omega]_i$ times in non-circled cells with one
additional time in a circled cell if $C[\Omega]_i$ is fermionic.  We
will call the set of such fillings $\mathcal
T^{(e)}[\Omega;\Lambda']$.

Finally, given a filling $T \in \mathcal T^{(e)}[\Omega;\Lambda']$,
we read the content of the circles from top to bottom and obtain a
word $a_1\dots a_m$. The sign of the permutation needed to reorder
this word such that it be {\it increasing} gives the weight
associated to the filling $T$, denoted this time $\bar{w}[T]$.  The
weight of $T$ is the sign needed to reorder the monomial associated
to $T$ so that it corresponds
  to   $\theta_{C[\Omega]} x^{C[\Omega]}$ up to a factor $(-1)^{m(m-1)/2}$.
This is because to coincide with the product in $e_{C[\Lambda]}$
being done columnwise, we would have to read from bottom to top.
Reading from top to bottom provides the $(-1)^{m(m-1)/2}$ factor
needed to obtain the coefficient in $\versg{e_\La}$ instead of in
$e_{C[\La]}$. We thus have \begin{equation}N^{\Om}_{\La}=\sum_{T \in
\mathcal T^{(e)}[\Omega;\Lambda']} \bar{w}[T]\, .\end{equation} We
now use this equation to prove the theorem.

First, it is easy to convince ourselves that there is only one
element in $\mathcal T^{(e)}[\Lambda';\Lambda']$ and that it has a
positive weight. Because the rows of $D[\La']$ and $C[\La']$
coincide, for the filling to have increasing rows, we have no choice
but to put the $C[\La']_i$ letters $i$ in the $i$-th row of
$D[\La']$.  In the case when $C[\La']_i$ is fermionic, the extra $i$
has no choice but to go in the circle in row $i$ of $D[\La']$ for no
two $i$'s to be in a same column.  For example, given
$\Lambda'=(3,1;2,1)$ filling $D[\La']$ with the letters of $C[\La']$
leads to:
\begin{equation}{ \tableau[scY]{1&1&1&\bl\tcercle{1}\\2&2\\3&
\bl\tcercle{3}\\4}}
\quad\mbox{in}\quad{\tableau[scY]{&&&\bl\tcercle{}\\&\\&
\bl\tcercle{}\\&\bl}}\quad \longrightarrow\quad
{\tableau[scY]{1&1&1&\bl\tcercle{1}\\2&2\\3& \bl\tcercle{3}\\4}}
\end{equation}

This explains the first term in (\ref{einm}).

Second, let $\om=\Om^*$  and $\la=\La^*$. If $\Om \not <_S\La'$, a
filling of $\Om$ by $\La$ is obviously impossible because we would
need to be able to obtain in particular (forgetting about the
circles) a filling of the type $\mathcal T^{(e)}[\om;\la']$, which
would contradict the well known fact that the theorem holds in the
zero-fermion case.

Finally, we suppose that $\La^*=\Om^*$ and ${\rm sh}(D[\Om]) \not <
{\rm sh}(D[\La'])$ (see Corollary~\ref{coropart}). Again, this is
impossible because we would need to be able to obtain in particular
(considering the circles as usual cells and reordering the columns)
a filling of the type $\mathcal T^{(e)}[{\rm sh}(D[\Om]);{\rm
sh}(D[\La'])]$, which would contradict the well known fact that the
theorem holds in the zero-fermion case.
\end{proof}
\n Note that for the various examples that we have worked out, the
coefficients  $N^{\Om}_{\La}$ happened to be non-negative.  So we
may surmise that a stronger version of the theorem, where
$N^{\Om}_{\La}$ is a non-negative integer, holds.

    The linear independence of the $e_\La$'s in
$\mathscr{P}^{S_N}$ implies that the first $N$ bosonic and fermionic
elementary  functions are algebraically independent over
$\mathbb{Z}$. Symbolically, \begin{equation}\label{fondtheo}
\mathbb{Z}[x_1,\dots,x_N,\theta_1,\ldots,\theta_N]^{S_N}\equiv
\mathbb{Z}[e_1,\dots,e_N,\tilde{e}_0,\dots\tilde{e}_{N-1}]\,
,\end{equation} which can be interpreted as the fundamental theorem
of symmetric  polynomials in superspace.

\subsection{Complete symmetric functions and involution}

The $n$-th {\it bosonic and fermionic  complete symmetric
 functions} are given respectively by
\begin{equation}\label{bhn} h_n:=\sum_{\lambda \vdash n}\,m_\lambda
\qquad {\rm and}\qquad \tilde{h}_{n}:=\sum_{\Lambda
\vdash(n|1)}(\La_1+1)\,m_\La\,,\end{equation}
   From the explicit form of $h_n(x)$, namely,
$\sum_{1\leq i_1 \leq i_2 \leq \ldots \leq i_n}x_{i_1}\cdots
x_{i_n}$, we see that its fermionic partner is again generated by
the action of $d$ in the form-representation:
\begin{equation}\tilde{h}_{n-1}(x,\theta)\sim
\tilde{h}_{n-1}(x,dx)=d\,h_n(x)\quad\mbox{for all}\quad n\geq 1\,
.\end{equation}

The generating function for complete symmetric  polynomials is
\begin{equation}\label{defH} H(t,\tau):=\sum_{n=0}^\infty
t^n(h_n +\tau\tilde{h}_n )=\prod_{i=1}^\infty \frac{1}{1-t x_i -
\tau \theta_i}\,.
\end{equation}
To prove (\ref{defH}), one simply uses the inversion of even
elements in the Grassmann algebra, which gives
  \begin{equation}\frac{1}{1-t x_i
- \tau \theta_i}= \sum_{n\geq 0}\left(t x_i+\tau \theta_i\right)^n=
\sum_{n\geq 0}\left[(t x_i)^n+n \tau \theta_i (tx_i)^{n-1}\right]\,
.\end{equation}

   From relations (\ref{defE}) and (\ref{defH}), we get
\begin{equation}H(t,\tau)\,E(-t,-\tau)=1\,.\end{equation} By
expanding the generating functions in terms of  $e_n$,
$\tilde{e}_n$, $h_n$ and $\tilde{h}_n$ in the last equation, we
obtain recursion relations, of which the non-fermionic one is a well
known formula.

\begin{lemma}\label{recureh} Let $n\geq 1$, then
\begin{equation} \sum_{r=0}^n(-1)^r e_r h_{n-r}=0\, .\end{equation}
Let $n \geq 0$, then\begin{equation} \sum_{r=0}^n(-1)^r
(e_r\tilde{h}_{n-r}-\tilde{e}_r h_{n-r})=0\,
.\end{equation}\end{lemma}

\n Note that the second relation can be obtained from the first one
by the action of $d$ (representing, as usual, $\ta_i$ as $dx_i$).

We consider a homomorphism  $\hat{\omega}\,
:\,\mathscr{P}^{S_\infty}(\mathbb{Z})\rightarrow
\mathscr{P}^{S_\infty}(\mathbb{Z})$ defined by the following
relations:
    \begin{equation}\label{definvolution}
    \hat{\omega} \,:\,
e_n\longmapsto h_n \quad\mbox{and}\quad
\tilde{e}_n\longmapsto\tilde{h}_n\, .\end{equation}
\begin{theorem}\label{involution}The homomorphism  $\hat{\omega}$ is
an involution, i.e., $\hat{\omega}^2=1$.  Equivalently, we have
\begin{equation}  \hat{\omega} \,:\, h_n\longmapsto e_n \quad\mbox{and}\quad
\tilde{h}_n\longmapsto\tilde{e}_n\,.\end{equation}\end{theorem}
\begin{proof} This comes from the application of transformation
(\ref{definvolution}) to the recursions appearing in
Lemma~\ref{recureh} followed by the comparison with the original
recursions.  Explicitly:
    \begin{equation}0=\sum_{r=0}^n(-1)^r
\hat{\omega}(e_r) \hat{\omega}(h_{n-r})=(-1)^n\sum_{r=0}^n(-1)^r
\hat{\omega}(h_{r})h_{n-r}\, ,\end{equation} which implies
$\hat{\omega}(h_{r})=e_r$. Similarly, we have
\begin{equation}0=\sum_{r=0}^n(-1)^r
\left(\hat{\omega}(e_r)\hat{\omega}(\tilde{h}_{n-r})-\hat{\omega}(\tilde{e}_r)
\hat{\omega}(h_{n-r})\right)=(-1)^{n-1}\sum_{r=0}^n(-1)^r \left(
e_r\tilde{h}_{n-r}-\hat{\omega}(\tilde{h}_{r})h_{n-r}\right)\,
,\end{equation} leading to
$\hat{\omega}(\tilde{h}_{r})=\tilde{e}_{r}$.\end{proof}

Now, let
\begin{equation}h_\La:=\prod_{i=1}^{\underline{\overline{\La}}}\tilde{h}_{\La_{i}}\prod_{j=\underline{\overline{\La}}+1}^{\ell(\La)}
{h}_{\La_{j}}\, . \end{equation} Equation (\ref{definvolution}) and
Theorem \ref{involution} immediately give a bijection between two
sets of multiplicative  polynomials:
\begin{equation}\hat{\omega}(e_\La)=h_\La\quad\mbox{and}\quad\hat{\omega}(h_\La)=e_\La\,
. \end{equation} We have thus obtained another $\mathbb{Z}$-basis
for the algebra of symmetric superpolynomials.
\begin{corollary} The set $\{\,h_\La \,:\, \La\vdash(n|m)\}$ is a basis of
$\mathscr{P}^{S_\infty}_{(n|m)}(\mathbb{Z})$.\end{corollary}

Finally, Lemma \ref{recureh} allows us to write determinantal
expressions for the elementary symmetric  functions in terms of the
complete symmetric  functions and {\it vice versa} using the
homomorphism $\hat \omega$.
\begin{proposition} \label{deteh} For $n\geq 1$, we have
\begin{equation}e_n=\left|%
\begin{array}{cccccc}
{h}_1  & {h}_2 & {h}_{3}&\cdots&h_{n-1}&{h}_n  \\
    1 & h_1 & h_2 & \cdots & h_{n-2}&h_{n-1} \\
    0 & 1 & h_1 & \cdots &h_{n-3}&h_{n-2}\\
    0&0&1&\ddots&h_{n-4}&h_{n-3}\\
     \vdots & \vdots & \ddots & \ddots& \ddots&\vdots \\
     0 & 0 & 0 & \cdots &1&h_1 \\
\end{array}%
\right|.\end{equation}For $n\geq 0$, we have
\begin{equation}\label{deuxiemedet}
\tilde{e}_n=\frac{1}{n!}\left|%
\begin{array}{cccccc}
     \tilde{h}_0 & \tilde{h}_1  & \tilde{h}_2 & \cdots
&\tilde{h}_{n-1}&\tilde{h}_n  \\
    n & (n+1)h_1 & (n+2)h_2 & \cdots & (2n-1) h_{n-1}&2n h_{n} \\
    0 & n-1 & n h_1 & \cdots &(2n-3) h_{n-2} &(2n-2)h_{n-1} \\
    0&0&n-2&\ddots&(2n-5) h_{n-3}&(2n-4)h_{n-2} \\
     \vdots & \vdots & \ddots & \ddots&\ddots&\vdots \\
     0 & 0 & 0 & \cdots & 1&2h_1 \\
\end{array}%
\right|.\end{equation}\end{proposition} \begin{proof} The first
relation is well known to be a simple application of Cramer's rule
to the linear system coming from Lemma~\ref{recureh}: $\mathbf{h}
=\mathbf{e}  \, \mathbf{H}$, where
\begin{equation} \mathbf{h}^\mathrm{t}=
\left(\begin{array}{c}
   h_1 \\ h_2 \\h_3\\ \vdots \\ h_N \\
\end{array}\right)\quad
\mathbf{e}^\mathrm{t}=\left(\begin{array}{c}
   e_1 \\ e_2 \\e_3\\ \vdots \\ e_N \\
\end{array}\right)
\quad \mathbf{H}=\left(\begin{array}{rrrrr}
   1&h_1&h_2&h_3&\ldots \\
   0&-1&-h_1&-h_2&\ldots  \\
   0&0&1&h_1&\ddots \\
   0&0&0&-1&\ddots\\
   \vdots & \ddots& \ddots& \ddots&\ddots\\
\end{array}\right).\end{equation}

To obtain the other determinant, we use the second formula of
Lemma~\ref{recureh} to obtain the linear system:
\begin{equation}
\tilde{\mathbf{h}} \, \overline{\mathbf{H}}=\tilde{\mathbf{e} }\,
\mathbf{H}\, ,\end{equation} where $\mathbf{H}$ is as given above,
and where
\begin{equation} \tilde{\mathbf{h}}^\mathrm{t}=
\left(\begin{array}{c}
   \tilde h_0 \\ \tilde h_1 \\ \tilde h_2\\ \vdots \\ \tilde h_{N-1} \\
\end{array}\right)\quad
\tilde{\mathbf{e}}^\mathrm{t} =\left(\begin{array}{c}
   \tilde e_0 \\ \tilde e_1 \\ \tilde e_2\\ \vdots \\ \tilde e_{N-1} \\
\end{array}\right)  \quad
  \overline{\mathbf{H}}=\left(\begin{array}{rrrrr}
   1&-e_1&e_2&-e_3&\ldots \\
   0&1&-e_1&e_2&\ldots  \\
   0&0&1&-e_1&\ddots \\
   0&0&0&1&\ddots\\
   \vdots & \ddots& \ddots& \ddots&\ddots\\
\end{array}\right). \end{equation}
Using the coadjoint formula for the inverse of a matrix, and the
determinantal expression for $h_n$ obtained by applying the
homomorphism $\hat \omega$ on the determinant of $e_n$ given above,
it is not hard to see that the $(i,j)$-th component of the inverse
of $\overline{\mathbf{H}}$ is simply $h_{j-i}$.  We are thus led to
the matrix relation:
\begin{equation}\tilde{\mathbf{h}}=\tilde{\mathbf{e}} \, \widetilde{\mathbf{H}} \,
,\end{equation} where
\begin{equation}\widetilde{\mathbf{H}}_{i,j}=(-1)^{i+1}\sum_{k\geq0}^{i-j}h_{i-j-k}h_k=
-\widetilde{\mathbf{H}}_{i+1,j+1} \, .\end{equation} If we set
$H_i=\sum_{k=0}^i h_{i-k} h_k$, the matrix $\widetilde{\mathbf{H}}$
  can be expressed in a convenient form as
\begin{equation} \widetilde{\mathbf{H}}=\left(\begin{array}{rrrrr}
   1&H_1&H_2&H_3&\ldots\\
   0&-1&-H_1&-H_2&\ldots\\
   0&0&1&H_1&\ldots \\
   0&0&0&-1&\ldots\\
  \vdots & \ddots& \ddots& \ddots&\ddots\\
\end{array}\right) .\end{equation}
Using Cramer's rule and then multiplying rows $2,3,\dots n$ of the
resulting determinant by $n,n-1,\dots,1$ respectively, we obtain
\begin{equation}\label{finaldet}
\tilde{e}_n=\frac{1}{n!}\left|%
\begin{array}{cccccc}
     \tilde{h}_0 & \tilde{h}_1  & \tilde{h}_2 & \cdots
&\tilde{h}_{n-1}&\tilde{h}_n  \\
    n & nH_1 & n H_2 & \cdots & n H_{n-1}&n  H_{n} \\
    0 & n-1 & (n-1) H_1 & \cdots &(n-1) H_{n-2} &(n-1)H_{n-1} \\
    0&0&n-2&\ddots&(n-2) H_{n-3}&(n-2)H_{n-2} \\
     \vdots & \vdots & \ddots & \ddots&\ddots&\vdots \\
     0 & 0 & 0 & \cdots & 1&H_1 \\
\end{array}%
\right|.\end{equation} We will finally show that by manipulating
determinant (\ref{deuxiemedet}), we obtain determinant
(\ref{finaldet}). Let $R_i$ be row $i$ of determinant
(\ref{deuxiemedet}). If $R_2 \to R_2+R_3 h_1+\dots+R_{n} h_{n-2}$ in
this determinant, the second row becomes that of determinant
(\ref{finaldet}) due to the simple identity (for $i=1,\dots, n$)
$$
nH_i = (n+i) h_i +(n+i-2) h_{i-1}h_1+ \cdots + (n+i-2i+2)h_1
h_{i-1}+(n+i-2i)h_i \, .
$$
Doing similar operations on the lower rows, the two determinants are
seen to coincide.
\end{proof}

\subsection{Power sums}
We define the $n$-th {\it bosonic and fermionic  power sums} as
follows: \begin{equation}p_n:=\sum_{i=1}^\infty x_{i}^n=m_{(n)}
\qquad {\rm and}\qquad \tilde{p}_n:=\sum_{i=1}^\infty\theta_i
x_i^n=m_{(n;0)}\, . \end{equation} Note that this time we will set
$p_0=0$. Obviously,
\begin{equation}n\,\tilde{p}_{n-1}(x,\theta)\sim
n\,\tilde{p}_{n-1}(x,dx)=d\,p_n(x)\end{equation}for all $n \geq 1$.

Proceeding as in the complete symmetric functions case, we introduce
products of power sums:
\begin{equation}
p_\La:=\prod_{i=1}^{\underline{\overline{\La}}}\tilde{p}_{\La_{i}}\prod_{j=\underline{\overline{\La}}+1}^{\ell(\La)}
{p}_{\La_{j}}\, .\end{equation} Also, we find that the generating
function for superpower sums is
\begin{equation}\label{defP}P(t,\tau):=\sum_{n\geq 0}t^n\,p_n
+\tau\sum_{n\geq 0}(n+1)t^n \,\tilde{p}_n =\sum_{i=1}^\infty \frac{t
x_i + \tau \theta_i}{1-t x_i - \tau \theta_i}\, .\end{equation} One
directly verifies
\begin{equation}H(t,\tau)\,P(t,\tau)=(t\partial_t+\tau\partial_\tau)H(t,\tau)\end{equation}
and\begin{equation}E(t,\tau)\,P(-t,-\tau)=-(t\partial_t+\tau\partial_\tau)E(t,\tau)\,
.\end{equation} These expressions  lead after some manipulations to
the following recursion relations.
\begin{lemma}\label{recurep} Let $n\geq 1$.  Then
\begin{equation} h_n=\sum_{r=1}^n p_r h_{n-r}\,
, \qquad  ne_n=\sum_{r=1}^n (-1)^{r+1}p_r e_{n-r}\, .
\end{equation}
Let $n\geq 0$, and recall that $p_0=0$.  Then
\begin{equation} (n+1)\,
\tilde{h}_n=\sum_{r=0}^n \bigl[ \, p_r\tilde{h}_{n-r}+
(r+1)\tilde{p}_r h_{n-r} \, \bigr]\, ,
\end{equation}
\begin{equation}(n+1)\,
\tilde{e}_n=\sum_{r=0}^n (-1)^{r+1}  \bigl[ \, p_r\tilde{e}_{n-r}-
(r+1)\tilde{p}_r e_{n-r} \, \bigr] \, .\end{equation}\end{lemma}

\begin{theorem}Let $\hat{\omega}$ be the involution defined in
(\ref{definvolution}). Then, for $n>0$,
    \begin{equation}\hat{\omega}\,:\,
p_n\longmapsto(-1)^{n-1}p_n\quad\mbox{and}\quad
     \tilde{p}_{n-1}\longmapsto(-1)^{n-1}\tilde{p}_{n-1}\end{equation}
or, equivalently,\begin{equation}\label{involup}
\hat{\omega}(p_\Lambda)=\omega_\La\,p_\Lambda\quad\mbox{with}\quad
\omega_\La:=(-1)^{|\La|+\overline{\underline{\Lambda}}-{\ell}(\La)}\,
.\end{equation}
\end{theorem}
\begin{proof}  We use  Lemma \ref{recurep} and  proceed as in the
proof of  Theorem \ref{involution}.\end{proof}

\begin{proposition} \label{detpe}For $n\geq 1$, we have
\begin{equation}p_n=\left|%
\begin{array}{ccccc}
{e}_1  & 2{e}_2 & 3{e}_{3}&\cdots&n{e}_n  \\
    1 & e_1 & e_2 & \cdots &e_{n-1} \\
    0 & 1 & e_1 & \cdots &e_{n-2}\\
     \vdots & \ddots & \ddots& \ddots&\vdots \\
     0 & 0 & 0&&e_1 \\
\end{array}%
\right|, \quad n!\,e_n=\left|%
\begin{array}{ccccc}
{p}_1  & p_2 &\cdots&p_{n-1}&p_n  \\
    1 & p_1  & \cdots &p_{n-2}&p_{n-1} \\
    0 & 2& \ddots &p_{n-3}&p_{n-2}\\
     \vdots & \ddots & \ddots& \ddots &\vdots \\
     0 & 0 &\cdots&n-1&p_1 \\
\end{array}%
\right|.\end{equation}For $n\geq 0$, we have
\begin{equation}\tilde{p}_n=\left|%
\begin{array}{ccccc}
     \tilde{e}_0 & \phantom{n}\tilde{e}_1  & \phantom{n}\tilde{e}_2 &
\cdots&\phantom{n}\tilde{e}_n  \\
    1 &\phantom{n} e_1 & \phantom{n}e_2 & \cdots &e_{n} \\
    0 & \phantom{n}1 & \phantom{n}e_1 & \cdots &e_{n-1} \\
     \vdots &  \ddots & \ddots&\ddots&\vdots \\
     0 & \phantom{n}0 & \phantom{n}0 & &e_1 \\
\end{array}%
\right|
,\quad {n!\,}\tilde{e}_n=\left|%
\begin{array}{ccccc}
     \tilde{p}_0 & \phantom{n}\tilde{p}_1   &
\cdots&\tilde{p}_{n-1}&\tilde{p}_n  \\
    n & \phantom{n}p_1 &  \cdots &p_{n-1}&p_{n} \\
    0 & n-1 &  \ddots &p_{n-2}&p_{n-1} \\
     \vdots &  \ddots & \ddots&\ddots &\vdots \\
     0 & \phantom{n}0 & \cdots & 1&p_1 \\
\end{array}%
\right|.\end{equation} Similar formulas for the complete symmetric
 functions are obtained by using the involution
$\hat\omega$.\end{proposition}
    \begin{proof} The proof is similar to that of Proposition
\ref{deteh}. \end{proof}

The  explicit formulas presented in Proposition \ref{detpe}
establish the correspondence between the sets
$\{p_n,\tilde{p}_{n-1}\}$ and $\{e_n,\tilde{e}_{n-1}\}$. This
implies, in particular, that $e_\La=\sum_\Om c_{\La\Om}\,p_\Om$ for
uniquely determined coefficients  $c_{\La\Om}\in\mathbb{Q}$. Note
that $c_{\La \Om}$ is not necessarily an integer since, for
instance, $e_2=p_1^2/2-p_2$. Theorem \ref{theoebase} and Proposition
\ref{detpe} thus imply the following result.

\begin{corollary}\label{theobasep} The set $\{\,p_\La :
\La\vdash(n|m)\}$ is a basis of
$\mathscr{P}^{S_N}_{(n|m)}(\mathbb{Q})$ if  either $n\leq N$ and
$m=0$ or $n<N$ and $m>0$.  In particular, the set $\{\,p_\La :
\La\vdash(n|m)\}$ is a basis of
$\mathscr{P}^{S_\infty}_{(n|m)}(\mathbb{Q})$.
\end{corollary}

 The
power sums will play a fundamental role in the remainder of the
article. For this reason, we will consider, from now on, only
symmetric  polynomials in superspace defined over the rational
numbers (or any greater field):
\begin{equation} \mathscr{P}^{S_\infty}:=  \mathscr{P}^{S_\infty}(\mathbb{Q})\,
.\end{equation}

\begin{remark}
The elementary and power-sum superpolynomials appear to have been
introduced first in \cite{brink}. The
 power-sum
superpolynomials were rediscovered in  \cite{DLM5} (section 2.5) and
the expressions for the three multiplicative bases were first given
in \cite{DLMpro}.
\end{remark}

\subsection{Orthogonality}
Let $n_\la(i)$ denote the number of parts  equal to $i$ in the
partition $\lambda$.  We introduce a bilinear form, $\LL\,|\,\RR :
\mathscr{P}^{S_\infty}\times
\mathscr{P}^{S_\infty}\rightarrow\mathbb{Q}$, defined by
\begin{equation}\LL \versg{p_\La} | \versd{p_\Om }\RR:=z_\La
\delta_{\La,\Om}\, , \quad {\rm for} \quad
z_\La:=z_{\La^s}=\prod_{k\geq
1}\left[\,k^{n_{\La^s}(k)}\,n_{\La^s}(k)!\,\right]\, .\end{equation}

\begin{proposition} Let $f$ and $g$ be superpolynomials in
$\mathscr{P}^{S_\infty}$.
    Then $\LL \versg{f}|\versd{g}\RR$ is a scalar product,
that is, in addition to bilinearity, we have
\begin{equation}\begin{array}{ll}  \LL
\versg{f}| \versd{g}\RR=\LL \versg{g}|\versd{f}\RR=\LL
\versd{g}|\versg{f}\RR&\mbox{(symmetry)}\\
\LL \versg{f}|\versd{f}\RR\,> \,0 \quad \forall\quad f\neq
0&\mbox{(positivity)}\, .\\
    \end{array} \end{equation}\end{proposition}
\begin{proof} The symmetry property is a consequence of Lemma
\ref{ordredestheta}.  The positivity of the scalar product is proved
as follows.  By definition  $z_\La>0$ and by virtue of Corollary
\ref{theobasep}, there is a unique decomposition $f=\sum_\La f_\La
p_\La$.  Therefore $\LL \versg{f} | \versd{f}\RR=\sum_\La f_\La^2
z_\La>0$. \hfill \end{proof}

\begin{proposition}The involution $\hat{\omega}$ is an
isometry.\end{proposition}
    \begin{proof} Given that
$\{p_\La\}_\La$ is a basis of $\mathscr{P}^{S_\infty}$, for any
symmetric polynomials $f$ and $g$, we have $f=\sum_\La f_\La p_\La$
and $g=\sum_\La g_\La p_\La$. Thus

\begin{eqnarray}
\LL \versg{\hat{\omega}f} |
\versd{\hat{\omega}g}\RR&=&\sum_{\La,\Om}f_\La g_\Om \LL
\versg{\hat{\omega}p_\La} | \versd{\hat{\omega}p_\Om}\RR\cr &=&\ds
\sum_{\La} f_\La g_\Om
(-1)^{|\La|+\overline{\underline{\Lambda}}-{\ell}(\La)}
(-1)^{|\Om|+\overline{\underline{\Om}}-{\ell}(\Om)}\LL\versg{p_\La}|\versd{p_\Om}\RR
\cr &=&\sum_{\La} z_\La f_\La g_\La=\LL \versg{f} | \versd{g}\RR\, ,
\end{eqnarray}as claimed.\end{proof}

The following theorem is of particular importance since it gives
Cauchy-type formulas for  the superpower sums.

\begin{theorem} \label{cauchypp} Let
$K=K(x,\theta;y,\phi)$   be the bi-symmetric formal  function given
by\begin{equation}\label{defK} K:=\prod_{
i,j}\frac{1}{1-x_iy_j-\theta_i\phi_j}  \; . \end{equation}
    Then
\begin{equation} K=\sum_{\La\in\mathrm{SPar}}z_\La^{-1}
\versg{p_\La(x,\theta)}\versd{p_\La(y,\phi)} \;.\end{equation}
\end{theorem}
    \begin{proof} We have:
\begin{eqnarray}
\prod_{i,j}\frac{1}{1-x_iy_j-\theta_i\phi_j}&=&\exp\Big\{\sum_{i,j}
\ln\Big[(1-x_iy_j-\theta_i\phi_j)^{-1}\Big]\Big\}
    =\exp\Big\{\sum_{i,j}\sum_{n\geq
    1}\Big[\frac{1}{n}(x_iy_j+\theta_i\phi_j)^{n}\Big]\Big\}\cr
    &=&\exp\Big\{\sum_{n\geq
1}\Big[\frac{1}{n}p_n(x)\,p_n(y)\Big]+\sum_{n\geq0}\Big[\tilde{p}_n(x,\theta)\,\tilde{p}_n(y,\phi)\Big]\Big\}\cr
    &=&\prod_{n\geq 1}\sum_{k_n\geq
    0}\frac{1}{n^{k_n}k_n!}\Big[p_n(x)\,p_n(y)\Big]^{k_n}\exp
\Big[\sum_{n\geq0}\tilde{p}_n(x,\theta)\,\tilde{p}_n(y,\phi)\Big]\,
    .
\end{eqnarray}
Considering Proposition \ref{expfermions}, we find
\begin{equation}\prod_{i,j}\frac{1}{1-x_iy_j-\theta_i\phi_j}=\sum_{n,\,m\geq
0}\sum_{\substack{\lambda\in \mathrm{Part}_s(n)\\ \mu\in
\mathrm{Part}_a(m)}} \Big[z_{\lambda}^{-1}p_\lambda(x)
p_\lambda(y)\versg{ \tilde{p}_\mu(x,\theta) }\versd{ \tilde{p}_\mu
(y,\phi) }\Big]\, .\end{equation} This equation (together with Lemma
\ref{ordredestheta}) proves the theorem. \end{proof}

\begin{remark}  The inverse of the kernel satisfies:
\begin{equation}
K(-x,-\theta;y,\phi)^{-1} =\prod_{i,j}(1+x_iy_j+\theta_i\phi_j) \,
=\sum_{\La\in\mathrm{SPar}}\omega_\La z_\La^{-1}
\versg{p_\La(x,\theta)}\versd{p_\La(y,\phi)}\, .
\end{equation}
The proof of this  result is similar to that of  Theorem
\ref{cauchypp}, apart from the presence of the coefficient
$\omega_\La=(-1)^{|\La|+\overline{\underline{\Lambda}}-{\ell}(\La)}$,
which comes from the expansion of $\ln(1+x_iy_j+\theta_i\phi_j)$.
This shows that
\begin{equation}\label{deuxK}
K(x,\theta;y,\phi) = \omega^{(x,\ta)}\, K(-x,-\theta;y,\phi)^{-1} =
\omega^{(y,\phi)}\, K(-x,-\theta;y,\phi)^{-1}
\end{equation}
where $\omega^{(x,\ta)}$ indicates that $\omega$ acts on the
$(x,\ta)$ variables and similarly for $\omega^{(y,\phi)}$.
\end{remark}

We now give two direct consequences of Theorem \ref{cauchypp}.
\begin{corollary} $K$ is a reproducing kernel
in the space of symmetric superfunctions:
\begin{equation}\LL\,K(x,\theta;y,\phi)\,|\,f(x,\ta\,\RR=f(y,\phi)\,
,\quad\mbox{for all}\quad f\in\mathscr{P}^{S_{\infty}}\,
.\end{equation}\end{corollary}

\begin{proof} If $f\in\mathscr{P}^{S_{\infty}}$, there exist unique
coefficients $f_\La$ such that $f=\sum_\La f_\La p_\La$. Hence,
\begin{eqnarray}
\LL\,K(x,\theta;y,\phi)\,|\,f(x,\ta)\,\RR
&=&\sum_{\Om,\La}z_\Om^{-1}f_\La\LL
\versg{p_\Om(x,\ta)}|\versd{p_\La(x,\ta)}\RR
\versd{p_\Om(y,\phi)}\cr &=&\sum_\La f_\La \versd{p_\La(y,\phi)}=
f(y,\phi)\, ,\end{eqnarray} as desired.\end{proof}

\begin{corollary} \label{einbasep}We have
\begin{eqnarray} h_n=\sum_{\Lambda\in \footnotesize{\spar}
(n|0)}z_\La^{-1}p_\Lambda\aand
    e_n=\sum_{\Lambda\in \footnotesize{\spar} (n|0)}z_\La^{-1}\omega_\Lambda
    p_\Lambda\, ,\end{eqnarray}
\begin{equation} \tilde{h}_n=\sum_{\Lambda\in\footnotesize{ \spar}
(n|1)}z_\La^{-1}p_\Lambda\aand
    \tilde{e}_n=\sum_{\Lambda\in \footnotesize{\spar}
(n|1)}z_\La^{-1}\omega_\Lambda
    p_\Lambda\, .\end{equation}\end{corollary}
\begin{proof} Using the definition of the generating function
$E(t,\tau)$, we  first make the following correspondence:
\begin{equation}E(t,0)=\sum_{n\geq 0}t^ne_n(x)=K(-x,0;y,0)^{-1}\Big|_{y=(t,0,0,\ldots)}\;.
\end{equation}
Thus, from Theorem \ref{cauchypp} and
$p_\la(y)|_{y=(t,0,0,\ldots)}=t^{|\lambda|}$, we have
\begin{equation}\sum_{n\geq0}t^ne_n(x)=\sum_{\lambda\in\mathrm{Par}}t^{|\lambda|}\omega_\lambda
z_\lambda^{-1} p_\lambda(x)\quad\Longrightarrow\quad
e_n=\sum_{\lambda\vdash n}z_\la^{-1}\omega_\lambda
    p_\lambda\, .\end{equation} Then, we observe that
    \begin{equation}\partial_\tau E(t,\tau)=\sum_{n=0}^{N-1}t^n
    \tilde{e}_n(x,\tau)=\partial_\tau
    K(x,\theta;y,\phi)\Big|_{\substack{y=(t,0,0,\ldots)\\
    \phi=(-\tau,0,\ldots)}}\; .\end{equation} Theorem
    \ref{cauchypp} and
\begin{equation}
p_\La(y,\phi)\Big|_{\substack{y=(t,0,0,\ldots)\\
    \phi=(-\tau,0,\ldots)}}= \left \{
\begin{array}{ll} \phantom{-\tau}t^{|\La|} & {\rm ~if~} \overline{\underline{\La}}=0 \\
-\tau t^{|\La|} & {\rm ~if~} \overline{\underline{\La}}=1  \\
   \phantom{-\tau}0 & {\rm ~otherwise~}
\end{array} \right.\phantom{\{}
\end{equation}
finally lead to
    \begin{equation}\sum_{n\geq0}t^n
\tilde{e}_n(x,\tau)=\sum_{\Lambda,\overline{\underline{\La}}=1}z_\La^{-1}\omega_\Lambda
    p_\Lambda(x,\theta)\quad
    \Longrightarrow\quad \tilde{e}_n=\sum_{\Lambda\in
\mathrm{SPar}(n|1)}z_\La^{-1}\omega_\Lambda
    p_\Lambda\, .\end{equation}
    Note that the minus sign disappears since $\partial_\tau$ and
$p_\La(x,\theta)$ anticommute when $\overline{\underline{\La}}=1$.
    Similar formulas relating the superpower sums to the
homogeneous symmetric  polynomials
    are obtained using the involution $\hat{\omega}$.\end{proof}

\begin{lemma}\label{cauchybases}Let $\{u_\La\}$ and $\{v_\La\}$ be two
bases of $\mathscr{P}^{S_{\infty}}_{(n|m)}$.  Then
\begin{equation}K(x,\theta;y,\phi)=\sum_\La
\versg{u_\La(x,\theta)}\,\versd{v_\La(y,\phi)}\quad\Longleftrightarrow\quad\LL
\versg{u_\La}|\versd{v_\La}\RR=\delta_{\La,\Om}\,
.\end{equation}\end{lemma}
\begin{proof}
The proof is identical to the one in the case without Grassmannian
variables (see \cite{Mac} I.4.6).\end{proof}
% Given that $\{u_\La\}$ and $\{v_\La\}$ are bases,
% we can write
% \begin{equation}\versg{p_\La}=z_\La\sum_{\Om}U_{\La\Om}\versg{u_\Om}
% \aand \versd{p_\La}=\sum_{\Om}V_{\La\Om}\versd{v_\Om}\end{equation}
% This implies, on the one
% hand,\begin{equation}\delta_{\La,\Om}=z_\La^{-1}\LL
% \versg{p_\La}|\versd{p_\Om}\RR=
% \sum_{\Gamma,\Delta}U_{\La\Gamma}V_{\Omega\Delta}\LL
% \versg{u_\Gamma}|\versd{v_\Delta}\RR\end{equation} and
% \begin{equation}\delta_{\La,\Om}=\LL\versg{u_\La}|\versd{v_\Om}\RR
% \quad\Longleftrightarrow\quad\mathbf{1}=\mathbf{U}
% \mathbf{V}^{\mathrm{t}}\, ,\end{equation} where $\mathbf{U}$ is the
% matrix with entries $U_{\La\Om}$ and correspondingly for
% $\mathbf{V}$.  On the other hand, we have
% \begin{equation}K_+=\sum_\La z_\La^{-1}\versg{p_\La(x,\theta)}\,\versd{p_\La(y,\phi)}
% %
% =\sum_{\La,\Om,\Gamma}U_{\La\Om}V_{\La\Gamma}\vers% %      g{u_\Om(x,\theta)}\,\versd{v_\Gamma(y,\phi)}\,
% ,\end{equation} which means
% \begin{equation}K_+=\sum_\La
% \versg{u_\La(x,\theta)}\,\versd{v_\La(y,\phi)}
% \quad\Longleftrightarrow\quad\mathbf{1}=\mathbf{U}^{\mathrm{t}}\mathbf{V}\,
% .\end{equation} This completes the proof. \end{proof}

\begin{proposition} \label{propkernelmh} Let $K$ be the
 function
defined in (\ref{defK}). Then,
\begin{equation}K=\sum_{\La\in\mathrm{SPar}}
\versg{m_\La(x,\theta)}\versd{h_\La(y,\phi)} \, .\end{equation}
\end{proposition}

\begin{proof} We start with the definition of the generating
function $E(t,\tau)$:
\begin{eqnarray}
K(-x,-\ta;y,\phi)^{-1}&=&\prod_{ i,j}(1+x_iy_j+\theta_i\phi_j) =
\prod_{i }E(x_i, \ta_i) \cr &=&  \prod_{i
}\Big[\sum_{n\geq0}x_i^n\,e_n(y) +\theta_i\sum_{n\geq0}
x_i^n\tilde{e}_n(y,\phi)\Big]\cr
&=&\sum_{\epsilon_1,\epsilon_2,\dots \in \{0,1\}}
\sum_{n_1,n_2,\dots \geq 0} \prod_i \left(\theta_i^{\epsilon_i}
x_i^{n_i} e_{n_i}^{(\epsilon_i)}(y,\phi)  \right) \cr &=&
\sum_{\La\in\mathrm{SPar}}\versg{m_\La(x,\theta)}\versd{e_\La(y,\phi)}\;.
    \end{eqnarray}
In the third line we have set $e_n^{(0)}(y,\phi)=e_n(y,\phi)$ and
$e_n^{(1)}(y,\phi)=\tilde e_n(y,\phi)$. The fourth line follows by
reordering the variables using Lemma \ref{ordredestheta}. Using
(\ref{deuxK}), we can recover $K(x,\ta;y,\phi)$ by acting with
$\hat\omega^{\{y,\phi \}}$ on $ K(-x,-\ta;y,\phi)^{-1}$. The
identity then follows from $\hat\omega(e_{\Lambda})= h_{\Lambda}$.
\end{proof}

The previous proposition and Lemma \ref{cauchybases} have the
following corollary.
\begin{corollary}
The  monomials are dual to the complete symmetric
 functions in superspace:
\begin{equation} \LL \versg{\,h_\La} | \versd{m_\Om
}\RR=\delta_{\La,\Om}\, .\end{equation}\end{corollary}

\section{Concluding remarks}

% \subsection{Summary}

Most of the key classical concepts in the theory of symmetric
functions have been extended to superspace. The notable exception
concerns  the Schur function $s_\la$. The proper superspace
generalization of the classical definition of Schur functions as a
bialternant has not been found yet. (We point out in that regard
that division by anticommuting variables is prohibited). Similarly,
the Jacobi-Trudi identity, which expresses the Schur functions in
terms of the $h_n$'s, has not been generalized. Notice  that in all
the instances where we have obtained a determinantal expression, we
had at most one row or one column  made out of fermionic quantities,
something which cannot be the case for the sought for Jacobi-Trudi
super-identity. To determine whether these properties are specific
to the $m=0$ sector or not requires further study.

Note also that, off-hand, it appears unlikely that the
yet-to-be-defined Schur
 functions in superspace would be related to the representation theory of
special Lie superalgebras since these theories do not involve
Grassmannian variables. (Recall in that regard that the special
importance of the Schur functions lies in their deep representation
theoretic interpretation: $s_\la$ is a Lie-algebra character,  being
expressible as a sum of semistandard tableaux of shape $\la$.)
Actually, it could  well be that for the Schur superpolynomials, the
representation theoretic interpretation is simply lost.

The simplest way of defining the Schur superpolynomials is by a
specialization of the Jack superpolynomials. The latter have been
defined in \cite{DLM3}
 as the unique polynomials in $\mathscr{P}^{S_N} $ which are:  (1) eigenfunctions of supersymmetric differential operators; (2) triangular
 in the monomial basis with respect to the Bruhat order.
 It turns out that the Jack polynomials in superspace are also orthogonal  with respect to a one-parameter
 deformation of the  scalar product introduced in Section 3.4.  This will be considered
 in a forthcoming article \cite{DLM8b}.

In a different vein, with the introduction of diagrams with circles,
we expect a large number of results linked to ``Ferrer-diagram
combinatorics'' to have nontrivial extensions to the  supercase.
Pieces of ``supercombinatorics'' have already been presented at the
end of Section 2.4.

\begin{acknow}
We thank A. Joyal for pointing out the connection between ${\tilde
e}_{n-1}, \, {\tilde h}_{n-1}$ and $ n{\tilde p}_{n-1}$ and the
corresponding  one-forms  $de_n,\, dh_{n} $ and $ dp_{n}$. This work
was  supported by NSERC and FONDECYT (Fondo Nacional de Desarrollo
Cient\'{\i}fico y Tecnol\'ogico) grant \#1030114. P.D. is grateful
to the Fondation J.A.-Vincent for a student fellowship and to NSERC
for a postdoctoral  fellowship.
\end{acknow}

\end{document}